
\input amstex
\documentstyle{amsppt}
\NoRunningHeads
\magnification = \magstep 1
\pagewidth{32pc}
\pageheight{48pc}

\topmatter
\title\  The equality
$I^2=QI$ in Buchsbaum rings
\endtitle

\abstract Let $A$ be a Noetherian local
ring with the maximal ideal $\frak{m}$ and
$d = \roman{dim}~A$. Let $Q$ be a
parameter ideal in $A$. Let $I = Q
: \frak{m}$. The problem
of when the equality $I^2=QI$ holds true
is explored. When $A$ is a
Cohen-Macaulay ring, this problem was
completely solved by A. Corso, C.
Huneke, C. Polini, and W. Vasconcelos
\cite{CHV, CP, CPV}, while
nothing is known when
$A$ is not a Cohen-Macaulay ring. The
present purpose is to show that within a
huge class of Buchsbaum local rings $A$
the equality $I^2=QI$ holds true for all
parameter ideals $Q$. The result will
supply
\cite{Y1, Y2} and \cite{GN} with ample
examples of ideals $I$, for which the
Rees algebras $\roman{R}(I)= \bigoplus_{n
\geq 0}I^n$, the associated graded
rings $\roman{G}(I) =
\roman{R}(I)/I\roman{R}(I)$, and the
fiber cones
$\roman{F}(I)=\roman{R}(I)/\frak{m}\roman{R}(I)$
are all Buchsbaum rings with certain
specific graded local cohomology modules.
Two examples are explored. One is to
show that $I^2=QI$ may hold true for all
parameter ideals
$Q$ in $A$, even though $A$ is not a
generalized Cohen-Macaulay ring, and the
other one is to show that the equality
$I^2
= QI$ may fail to hold for some
parameter ideal
$Q$ in
$A$, even though $A$ is a Buchsbaum local
ring with multiplicity at least three.
\endabstract

\leftheadtext{Shiro Goto and Hideto
Sakurai}
\rightheadtext{The equality
$I^2=QI$ in Buchsbaum rings}

\author Shiro Goto and Hideto
Sakurai
\endauthor

\address Department of Mathematics, School of Science and
Technology, Meiji University, 214-8571 JAPAN 
\endaddress  

\email goto\@math.meiji.ac.jp \endemail

\address Department of Mathematics, School of Science and
Technology, Meiji University, 214-8571 JAPAN 
\endaddress  

\email ee78052\@math.meiji.ac.jp
\endemail

 

\thanks The first author is supported
by the Grant-in-Aid for Scientific
Researches in Japan (C(2), No. 13640044).
\endthanks
\endtopmatter

\def\depth{\operatorname{depth}}

\def\Ass{\operatorname{Ass}}

\def\Ker{\operatorname{Ker}}

\def\bar{\overline}


\def\z{{\Bbb Z}}

\def\z{{\Bbb Z}}

\font\b=cmr10 scaled \magstep4

\def\bigzerou{\smash{\lower1.7ex\hbox{\b 0}}}

\document
\baselineskip=15pt


\head 1. Introduction.
\endhead 

Let $A$ be a Noetherian local
ring with the maximal ideal $\frak{m}$
and $d =
\dim A$. Let $Q$ be a parameter ideal
in $A$ and let $I = Q : \frak{m}$. In
this paper we will study the problem of
when the equality $I^2 = QI$ holds
true. K. Yamagishi \cite{Y1, Y2} and
the first author and K. Nishida
\cite{GN} have recently showed the
Rees algebras
$\roman{R}(I)=
\bigoplus_{n
\geq 0}I^n$, the associated graded
rings $\roman{G}(I) =
\roman{R}(I)/I\roman{R}(I)$, and the
fiber cones
$\roman{F}(I)=\roman{R}(I)/\frak{m}\roman{R}(I)$
are all Buchsbaum rings with very specific
graded local cohomology modules, if
$I^2=QI$ and the base rings $A$ are
Buchsbaum. Our results will
supply 
\cite{Y1, Y2} and \cite{GN} with ample
examples.

Our research dates back to
the remarkable results of A. Corso, C.
Huneke, C. Polini, and W. Vasconcelos
\cite{CHV, CP, CPV}, who asserted that if
$A$ is a Cohen-Macaulay local ring,
then the equality
$I^2=QI$ holds true for every parameter
ideal $Q$ in $A$, unless $A$ is a regular
local ring. Let $\frak{a}^{\sharp}$ 
denote, for
an ideal
$\frak{a}$ in
$A$, the integral closure of $\frak{a}$.
Then their results are summarized into
the following, in which the equivalence
of assertions (2) and (3) are due to
\cite{G3, Theorem (3.1)}. The reader may
consult 
\cite{GH} for a simple proof of Theorem
(1.1) with a slightly
general form.

\proclaim{Theorem (1.1)
(\cite{CHV, CP, CPV})} Let
$A$ be a Cohen-Macaulay ring with
$\dim A = d$. Let $Q$ be a parameter
ideal in $A$ and let $I = Q : \frak{m}$.
Then the following three conditions are
equivalent to each other.

\roster
\item $I^2 \ne QI$.
\item $Q = Q^{\sharp}$.
\item $A$ is a regular local ring which
contains a regular system
$a_1, a_2,
\cdots, a_d$ of parameters such
that $Q = (a_1,\cdots, a_{d-1}, a_d^q)$
for some $1 \leq q \in \z$.
\endroster
Hence $I^2 = QI$ for every parameter
ideal $Q$ in $A$, unless $A$ is a regular
local ring.
\endproclaim

Our purpose is to
generalize Theorem (1.1) to local
rings $A$ which are not necessarily
Cohen-Macaulay. Since the notion of Buchsbaum ring
is a straightforward generalization of that of
Cohen-Macaulay ring, it seems quite
natural to expect that the equality
$I^2=QI$ still holds true also in
Buchsbaum rings. This is,
nevertheless, in general not true and
a counterexample is already explored by
\cite{CP}. Let $A = k[[X,Y]]/(X^2, XY)$
where
$k[[X,Y]]$ denotes the formal power
series ring in two variables over a field
$k$ and let $x, y$ be the reduction of
$X,Y$ mod the ideal $(X^2, XY)$. Let $Q =
(y^3)$ and put $I = Q : \frak{m}$. Then
$I = (x, y^2)$ and $I^2 \ne QI$ (\cite{CP,
p. 231}). However, the ideal $Q$ is
actually
$not$ the reduction of $I$ and
the multiplicity $\roman{e}(A)$ of
$A$ is 1. The Buchsbaum local ring
$A$ is $almost$ a DVR in the sense that
$A/(x)$ is a DVR and
$\frak{m}{\cdot}x=(0)$. Added
to it, with no difficulty one is able to check
that for a given parameter ideal $Q$ in
$A$, the equality $I^2= QI$ holds true if
and only if $Q
\not\subseteq \frak{m}^2$. For these
reasons this
example looks rather dissatisfactory, and
we shall provide in this paper more
drastic counterexamples. Nonetheless,
the example \cite{CP,
p. 231} was invaluable for the
authors to settle their starting point
towards the present research. For
instance, it strongly suggests that for
the study of the equality
$I^2=QI$ we first of all have to find the
conditions under which $Q$ is a
reduction of $I$, and the condition
$\roman{e}(A) \ne 1$ might play a certain
role in it. Any DVR contains
no parameter ideals
$Q$ for which the equality $I^2=QI$ holds
true, while as the example shows, non-Cohen-Macaulay Buchsbaum
local rings with
$\roman{e}(A) = 1$ could contain
somewhat ampler parameter ideals $Q$ for
which the equality $I^2=QI$ holds true.

Our problem is, therefore,
divided into two parts. One is to clarify
the condition under which $Q$ is
a reduction of $I$ and the other one is
to evaluate, when $I
\subseteq Q^{\sharp}$, the reduction number
$$\roman{r}_Q(I) = \min \{0 \leq n \in \z
\mid I^{n+1}= QI^n\}$$
of $I$ with respect to $Q$. As we shall quickly
show in this paper, one always has that
$I
\subseteq Q^{\sharp}$, unless
$\roman{e}(A) = 1$. In contrast, the
second part of our problem is in
general quite subtle and unclear, as we
will eventually show in this paper. We
shall, however, show that within a huge
class of Buchsbaum local rings $A$, the
equality
$I^2 = QI$ holds true for every
parameter ideal $Q$ in $A$.

Let us now state more
precisely our main results,
explaining how this paper is organized.
In Section 2 we will prove that if
$\roman{e}(A) > 1$, then $I= Q :
\frak{m}
\subseteq Q^{\sharp}$ for every parameter
ideal $Q$ in $A$. Hence $Q$ is a $minimal$
reduction of $I$, satisfying the
equality $\frak{m}I^n =
\frak{m}Q^n$ for all $ \in \z$
(Proposition (2.3)). Our proof is based
on the induction on $d = \dim A$, and the
difficulty that we meet
whenever we will check whether
$I^2 = QI$ is caused by the wild
behavior of the socle $(0) : \frak{m}$
in $A$. So, in Section 2, we shall
closely explain the method
how to control the socle $(0) :
\frak{m}$ in our context (Lemma (2.4)).
The main results of the section are Theorem
(2.1) and Corollary (2.13), which assert
that every unmixed local ring $A$ with
$\dim A \geq 2$ contains infinitely many
parameter ideals $Q$, for which the
equality $I^2 = QI$ holds true.  

In Section 3 we are concentrated to the
case where $A$ is a Buchsbaum local
ring. Let $A$ be a Buchsbaum local ring
with $d = \dim A \geq 1$ and let $x_1,
x_2, \cdots, x_d$ be a system of
parameters in $A$. Let $n_i \geq 1$ $(1
\leq i \leq d)$ be integers and put
$Q = (x_1^{n_1}, x_2^{n_2}, \cdots,
x_d^{n_d})$. We will then show that
$I^2=QI$ if $\roman{e}(A) >1$ and if $n_i
\geq 2$ for some $1 \leq i \leq d$
(Theorem (3.3)). Consequently, in a
Buchsbaum local ring $A$ of the form $A =
B/(f^n)$ where $n \geq 2$ and $f$ is a
parameter in a Buchsbaum local ring $B$,
the equality
$I^2=QI$ holds true for every parameter
ideal $Q$ (Corollary (3.7)). 

Let $\roman{r}(A) =
\underset{Q}\to{\sup}
\ \ \ell_A((Q :
\frak{m})/Q)$ where $Q$ runs over
parameter ideals in $A$, which we call
the Cohen-Macaulay type of $A$.
Then, thanks to Theorem (2.5) of
\cite{GSu}, one has the
equality
$$\roman{r}(A) = \sum_{i =
0}^{d-1}{d\choose i}h^i(A) +
\mu_{\hat{A}}(\roman{K}_{\hat{A}})$$
for every Buchsbaum local ring $A$ with
$d =
\dim A \geq 1$, where
$h^i(A) =
\ell_A(\roman{H}_{\frak{m}}^i(A))$
denotes the length of the
$i^{~\underline{th}}$ local cohomology
module of $A$ with respect to $\frak{m}$
and
$\mu_{\hat{A}}(\roman{K}_{\hat{A}})$
denotes the number of generators for the
canonical module
$\roman{K}_{\hat{A}}$ of the
$\frak{m}$-adic completion
$\widehat{A}$ of
$A$. Accordingly, one has $\ell_A((Q :
\frak{m})/Q) \leq
\roman{r}(A)$ in general, and if
furthermore $\ell_A((Q :
\frak{m})/Q) =
\roman{r}(A)$,
then the equality $I^2 = QI$ holds true
for the ideal $I = Q : \frak{m}$,
provided
$A$ is a Buchsbaum local ring  with
$\roman{e}(A) > 1$ (Theorem (3.9)).
Consequently, if $A$ is a Buchsbaum
local ring with $\roman{e}(A) >1$ and
the index $\ell_A((Q : 
\frak{m})/Q)$ of reducibility of $Q$ is
independent of the choice of a
parameter ideal
$Q$ in $A$, the equality $I^2=QI$ then
holds true for all parameter ideals $Q$
in
$A$. This result seems to account well
for the reason why Theorem (1.1) holds
true for Cohen-Macaulay rings $A$. In
Section 3 we shall also show that for a
Buchsbaum local ring $A$, there exists an
integer $\ell \gg 0$ such that the
equality $\roman{r}(A) = \ell_A((Q :
\frak{m})/Q)$ holds true for all
parameter ideals $Q \subseteq
\frak{m}^{\ell}$ (Theorem
(3.11)). Thus, inside
Buchsbaum local rings $A$ with $d = \dim
A \geq 2$, the parameter ideals $Q$ 
satisfying the equality $I^2 = QI$ are in
the majority. In the forthcoming paper
\cite{GSa} we will also prove that the
equality $I^2 = QI$ holds true for all
parameter ideals $Q$  in a Buchsbaum local
ring $A$ with $\roman{e}(A) = 2$
and $\depth A > 0 $.

In Section 4 we will give an
effective evaluation of the reduction
numbers
$\roman{r}_Q(I)$ in the case where $A$ is
a Buchsbaum local ring with $\dim A = 1$
and
$\roman{e}(A) > 1$ (Theorem (4.1)). The
evaluation is sharp, as we will show with
an example. The authors do not know
whether there exist some uniform
bounds of
$\roman{r}_Q(I)$ also in higher
dimensional cases.

It is somewhat
surprising to see that the equality
$I^2 = QI$ may hold true for $all$
parameter ideals $Q$ in $A$, even
though $A$ is not a generalized
Cohen-Macaulay ring. In Section 5 we
will explore one example satisfying
this property (Theorem (5.3)). In contrast, the
equality $I^2 = QI$ does in general not
hold true, even though $A$ is a
Buchsbaum local ring with $\roman{e}(A)
> 1$. In Section 5 we shall also
explore one more example of dimension 1
(Theorem (5.17)), giving complete
criteria of the equality $I^2 = QI$
for parameter ideals $Q$ in the
example.

We are now entering the very details.
Before that, let us fix again our
standard notation. Throughout, let
$(A,\frak{m})$ be a Noetherian local
ring with $d = \dim A$. We denote by
$\roman{e}(A) =
\roman{e}_{\frak{m}}^0(A)$ the
multiplicity of $A$ with respect to the
maximal ideal $\frak{m}$. Let
$\roman{H}_{\frak{m}}^i(*)$ denote the
local cohomology functor with respect
to $\frak{m}$. We denote by $\ell_A(*)$
and $\mu_A(*)$ the length and the number
of generators, respectively. Let
$\frak{a}^{\sharp}$ denote for an
ideal $\frak{a}$ in $A$ the
integral closure of $\frak{a}$. Let $Q
= (x_1, x_2, \cdots, x_d)$ be a
parameter ideal in $A$ and, otherwise
specified, we denote by $I$ the ideal 
$Q :
\frak{m}$. Let
$\roman{Min}~A$ be the set of minimal
prime ideals in
$A$. Let $\widehat{A}$ denote the
$\frak{m}$-adic completion of $A$.

\head 
2. A theorem for general local
rings.
\endhead

The goal of this section is the
following.

\proclaim{Theorem (2.1)}
Let $A$ be a Noetherian local ring with
$d = \dim A \geq 2$. Assume that $A$ is
a homomorphic image of a Gorenstein
local ring and $\dim A/\frak{p} = d$
for all $\frak{p} \in \roman{Ass}~A$.
Then $A$ contains a system $a_1, a_2,
\cdots, a_d$ of parameters such that for
all integers $n_i \geq 1$ $(1 \leq i
\leq d)$ the
equality $I^2 = QI$ holds true, where $$Q
=(a_1^{n_1}, a_2^{n_2}, \cdots,
a_d^{n_d}) \ \ \ \text{and} \ \ \ I = Q :
\frak{m}.$$ 
\endproclaim

To prove Theorem (2.1) we need some
preliminary steps. Let $A$
be a Noetherian local ring with the
maximal ideal $\frak{m}$ and
$d =\dim A \geq 0$. Let $Q$ be a
parameter ideal in $A$. We put $I = Q
:\frak{m}$. We begin with the
following.

\proclaim{Lemma (2.2)}
Suppose that $d \geq 1$. Then
$\roman{e}(A) = 1$ if $\frak{m}I
\not\subseteq \frak{m}Q$.
\endproclaim

\demo{Proof}
We may assume $I \ne A$. Let $W =
\roman{H}_{\frak{m}}^0(A)$ and $B =
A/W$. If $d = 1$, then $Q = \frak{m}I$,
since $Q$ is a principal ideal. Let $Q
= (a)$, $\bar{\frak{m}}= \frak{m}B$, and
$\bar{I} = IB$. Let $\bar{a} = a$ mod
$W$. Then, since
$(\bar{a})=\bar{\frak{m}}{\cdot}\bar{I}$
and $\bar{a}$ is a non-zerodivisor in
the Cohen-Macaulay local ring $B$, the
maximal ideal
$\bar{\frak{m}}$ is invertible, so that
$B$ is a DVR; hence $\roman{e}(B) =
\roman{e}(A) = 1$. Suppose that $d \geq
2$ and that our assertion holds true for
$d-1$. We choose $a_d \in \frak{m}I$ so
that $a_d
\not\in \frak{m}Q$, and then write $Q =
(a_1, \cdots, a_{d-1}, a_d)$. Let
$\bar{A} = A/(a_1)$, $\bar{\frak{m}} =
\frak{m}/(a_1)$,
$\bar{Q} = Q/(a_1)$, and $\bar{I} =
I/(a_1)$. Let $\bar{a_i}=a_i$ mod
$(a_1)$ $(2 \leq i \leq d)$. Then
$\bar{Q} = (\bar{a_2},
\cdots, \bar{a_d})$ is a parameter ideal
in $\bar{A}$ and
$\bar{I} = \bar{Q} : \bar{\frak{m}}$.
We have $\bar{\frak{m}}\bar{I}
\not\subseteq \bar{\frak{m}}\bar{Q}$,
since $\bar{a_d} \not\in
\bar{\frak{m}}\bar{Q}$. Hence
$\roman{e}(\bar{A}) = 1$ by the
hypothesis on $d$, so that $\roman{e}(A)
= 1$ as well.
\qed
\enddemo

\proclaim{Proposition (2.3)}
Suppose that $\roman{e}(A) > 1$. Then $I
\subseteq Q^{\sharp}$ and $\frak{m}I^n =
\frak{m}Q^n$ for all $n \in \z$.
\endproclaim

\demo{Proof}
We may assume that $d \geq 1$. Let $W =
\roman{H}_{\frak{m}}^0(A)$ and put $B =
A/W$. Then $\frak{m}B{\cdot}IB \subseteq
\frak{m}B{\cdot}QB$, since $\frak{m}I
\subseteq \frak{m}Q$ by Lemma (2.2).
Thus $IB$ is integral over
$QB$, because the ideal $\frak{m}B$
contains a non-zerodivisor of $B$
(recall that
$\depth B \geq 1$). Consequently, since
$W
\subseteq
\sqrt{(0)}$, $I$ is
integral over $Q$, so that $Q$ is a
minimal reduction of $I$. Since
$\frak{m}I \cap Q = \frak{m}Q$, we have
that $\frak{m}I = \frak{m}Q$, and
hence 
$\frak{m}I^n = \frak{m}Q^n$ for all $n
\in \z$.
\qed
\enddemo

The assertion that $I \subseteq Q^{\sharp}$ is
in general no longer true, unless
$\roman{e}(A) >1$ (see Theorem (1.1)).
When $A$ is not a Cohen-Macaulay ring,
the result is more complicated, as we
shall explore in Section 5.

The following result plays a key role throughout this paper as
well as in the proof of Theorem (2.1).

\proclaim{Lemma (2.4)}
Let $R$ be any commutative ring. Let
$M,L$, and $W$ be ideals in $R$ and let
$x
\in M$. Assume that $L : x^2 = L : x$
and $xW=(0)$. Then 
$$\left( L + (x^n) + W \right) : M =
\left[(L + W) : M \right] + \left[(L +
(x^n)) : M \right]$$
for all $n \geq 2$. If $L : x =
L : M$, we furthermore have that 
$$\left(L + (x^n) + W \right) : M
=\left(L+(x^n)\right) : M$$
for all  $n \geq 2$.
\endproclaim

\demo{Proof}
We have $L : x^{\ell} = L : x$ and
$[L+(x^{\ell})] \cap [L : (x^{\ell})] =
L$ for all $\ell \geq 1$, since $L : x^2
= L : x$. Let $\varphi \in (L + (x^n) +
W) : M$ and write $x\varphi = \ell + x^ny
+ w$, where $\ell \in L$, $y \in R$, and
$w \in W$. Let $z = \varphi
- x^{n-1}y$. Then since $x^2\varphi =
x\ell + x^{n+1}y$, we have $$z = \varphi
- x^{n-1}y \in L : x^2 = L : x. \tag 2.5$$
Let
$\alpha \in M$ and write $\alpha \varphi
= \ell_1 + x^ny_1 + w_1$ with $\ell_1 \in
L$, $y_1 \in R$, and
$w_1 \in W$. Then because 
$$\alpha \varphi = \ell_1 + x^ny_1 + w_1
= \alpha z + x^{n-1}(\alpha y)$$
we get $\alpha z - w_1 \in [L +
(x^{n-1})] \cap [L : x] \subseteq L$
(recall that $w_1 \in W \subseteq L : x$), whence
$$z \in (L + W) : M \subseteq (L + (x^n) +
W) : M$$ so that we also have
$x^{n-1}y =
\varphi - z \in (L + (x^n) + W) :
M$. Let
$\alpha \in M$ and write $x^{n-1}(\alpha
y) =\ell_2 + x^ny_2 + w_2$ with  $\ell_2
\in L$, $y_2 \in R$, and $w_2 \in W$.
Then $x^n(\alpha y) = x\ell_2 +
x^{n+1}y_2$ and $\alpha y -xy_2 \in L
: x^n = L : x$. Hence $y \in ((L : x) +
(x)) : M,$ so that $x^{n-1}y \in (L +
(x^n)) : M$ since $n \geq 2$. Thus 
$$\varphi = z + x^{n-1}y \in [(L + W) : M]
+[(L + (x^n)) : M].$$
If $L : x = L : M$ in addition, we get $z
\in L : M$ by (2.5), whence $$\varphi = z +
x^{n-1}y \in [L:M] + [(L + (x^n)):M] = (L
+ (x^n)) : M$$ as is claimed.
\qed
\enddemo

Let $R$ be a commutative ring and 
$x_1, x_2, \cdots, x_s \in R$ $(s \geq
1)$. Then $x_1, x_2, \cdots, x_s$ is
called a $d$-sequence in $R$, if $$(x_1,
\cdots, x_{i-1}) : x_i =(x_1,
\cdots, x_{i-1}) : x_ix_j$$ whenever $1
\leq i \leq j \leq s$. We say
that $x_1, x_2, 
\cdots, x_s$ forms a strong $d$-sequence
in $R$, if $x_1^{n_1},x_2^{n_2}, \cdots,
x_s^{n_s}$ is a $d$-sequence in $R$ for
all integers $n_i \geq 1~(1 \leq i \leq
s)$. See \cite{H} for basic but deep results on
$d$-sequences, which we shall freely use
in this paper. For example, if $x_1,
x_2, \cdots, x_s$ is a $d$-sequence in
$R$, then 
$$\align
(x_1, \cdots, x_{i-1}) : x_i^2
&=(x_1,
\cdots, x_{i-1}) : x_i  \tag
2.6\\
&=(x_1, \cdots,
x_{i-1}) : (x_1, x_2, \cdots, x_s) 
\endalign$$
for all $1 \leq i \leq s$.
Also one has the equality
$$\left((x_1, \cdots, x_{i-1}) :
x_i \right)
\cap (x_1, x_2, \cdots, x_s)^n = (x_1,
\cdots, x_{i-1}){\cdot}(x_1, x_2, \cdots,
x_s)^{n-1} \tag 2.7$$
for all integers $1 \leq i \leq s$ and $1
\leq n
\in \z$.

The following result is due to N. T.
Cuong.

\proclaim{Proposition (2.8) (\cite{C,
Theorem 2.6})} Let $A$ be a Noetherian
local ring with
$d = \dim A \geq 1$. Assume that $A$ is a
homomorphic image of a Gorenstein local
ring and that $\dim A/\frak{p} = d$ for
all $\frak{p} \in \Ass~A$. Then $A$
contains a system $x_1, x_2, \cdots, x_d$
of parameters which forms a strong
$d$-sequence.
\endproclaim

We will apply the following  result to 
strong $d$-sequences of Cuong.

\proclaim{Proposition (2.9)}
Let $R$ be a commutative ring and let
$x_1, x_2, \cdots, x_s \in R$ $(s \geq
1)$. Let
$Q = (x_1, x_2, \cdots, x_s)$ and $W
= (0) : Q$. Let $M$ be an ideal in $R$
such that $Q \subseteq M$. Assume that
$x_1, x_2, \cdots, x_s$ is a strong
$d$-sequence in $R$. Then
$$[(x_1^{n_1}, x_2^{n_2}, \cdots,
x_s^{n_s}) + W]: M = W + [(x_1^{n_1},
x_2^{n_2}, \cdots, x_s^{n_s}) : M]$$ for
all integers $n_i
\geq 2$
$(1 \leq i
\leq s).$
\endproclaim

\demo{Proof}
We put $L = (x_1^{n_1}, \cdots,
x_{s-1}^{n_{s-1}})$, $x = x_s$, and $n =
n_s$. Then $L : x^2 = L : x$, $x \in M$,
and $xW = (0)$. Hence by Lemma (2.4) we
get
$$[L + (x^n) + W] : M = [(L + W) : M] +
[(L + (x^n)) : M].\tag 2.10$$
Notice that $W : M = W$. (For, if
$\varphi \in W : M$, then $x_1 \varphi
\in W$ so that $x_1^2 \varphi = 0$,
whence
$\varphi \in (0) : x_1^2 = (0) : x_1 =
W$; cf (2.6).) Our
assertion is obviously true when $s =
1$.
Suppose that $s
\geq 2$ and that our assertion holds true
for $s-1$. Then, since $x_1, x_2, \cdots,
x_{s-1}$ is a strong $d$-sequence in $R$
and $W = (0) : x_1 = (0) : (x_1, \cdots,
x_{s-1})$ by (2.6), by the hypothesis
on $s$ we readily get that $$(L + W) : M =
W + (L : M)$$ whence by (2.10)
$$\align
[(x_1^{n_1}, x_2^{n_2}, \cdots,
x_s^{n_s}) + W] : M &= [(L + (x^n) + W)]
: M \\
&=[W + (L : M)] + [(L + (x^n)) : M]\\
&= W + [(L + (x^n)) : M]\\
&= W + [(x_1^{n_1}, x_2^{n_2}, \cdots,
x_s^{n_s}) : M]
\endalign$$
as is claimed.
\qed
\enddemo

We are now back to local rings.

\proclaim{Corollary (2.11)}
Let
$x_1, x_2,
\cdots, x_d$ be a
system of parameters in a Noetherian local ring $A$ with $d =
\dim A
\geq 1$
 and assume that $x_1, x_2,
\cdots, x_d$ forms a
strong
$d$-sequence. Let $n_i \geq 2$ $(1 \leq i \leq d)$
be integers and put $Q= (x_1^{n_1},
x_2^{n_2}, \cdots, x_d^{n_d}
)$. Then $I^2 = QI$ if
$\roman{e}(A) > 1$, where $I = Q :
\frak{m}.$  
\endproclaim

\demo{Proof}
Let $W = \roman{H}_{\frak{m}}^0(A)$.
Then $W = (0) : x_1 = (0) : (x_1, x_2,
\cdots, x_d)$. (For, if $\varphi \in W$,
then
$x_1^n \varphi = 0$ for some $n \gg 0$,
whence $\varphi \in (0) : x_1^n = (0) :
x_1 = (0) : (x_1, x_2, \cdots,
x_d)$; cf (2.6).) Let $B = A/W$. Then
since
$$(Q + W) : \frak{m} = W + (Q :
\frak{m}) = W + I$$
by
Proposition (2.9), we get $IB = QB :
\frak{m}B$. If $d = 1$, then $(IB)^2 =
QB{\cdot}IB$ by Theorem (1.1), because
$B$ is a Cohen-Macaulay ring with
$\roman{e}(B) =  \roman{e}(A) > 1$. Hence
$I^2 \subseteq QI + W$, so that we have 
$I^2 = QI$, because $$W \cap Q \subseteq
[(0) : (x_1)] \cap
(x_1, x_2, \cdots, x_d) = (0)$$ (cf.
(2.7)). Suppose that $d \geq 2$ and that
our assertion holds true for $d-1$. Let
$a_i = x_i^{n_i}$ $(1 \leq i \leq d)$
and put $\bar{A} = A/(a_1)$ and $\bar{I}
= I/(a_1)$. For each $c \in A$ let
$\bar{c}$ denote the reduction of $c$ mod
$(a_1)$. Then, since
$\roman{e}(\bar{A}) > 1$ and
the system $\bar{x_2},
\cdots, \bar{x_d}$ of parameters
for $\bar{A}$ forms by definition a
strong
$d$-sequence in $\bar{A}$, thanks to
the hypothesis on $d$, we get 
$\bar{I}^2 = (\bar{a_2}, \cdots,
\bar{a_d})\bar{I}$. Hence $I^2 \subseteq
(a_2, \cdots, a_d)I + (a_1)$ and so $I^2
= (a_2, \cdots, a_d)I + [(a_1) \cap I^2]$.

We then need the following.

\proclaim{Claim (2.12)}
$(a_1) \cap I^2 = a_1 I$.
\endproclaim

\demo{Proof of Claim (2.12)}
Let $\varphi \in (a_1) \cap I^2$ and
write $\varphi = a_1 y$ with $y \in A$.
Let $\alpha \in \frak{m}$. Then $\alpha
\varphi = a_1(\alpha y) \in Q^2$ since
$\frak{m} I^2 \subseteq Q^2$ (cf.
(2.3)). Consequently $a_1(\alpha y) \in
(a_1) \cap Q^2 = a_1Q$ (cf. (2.7)).
Hence
$\alpha y - q \in (0) : a_1 = (0) : x_1
= W$ for some $q \in Q$. Thus $$y \in (Q
+ W) : \frak{m} = W + I$$ so that
$\varphi = a_1y \in a_1 I$. Thus $(a_1)
\cap I^2 = a_1I,$ which completes the
proof of Corollary
(2.11) and  Claim (2.12) as well.
\qed
\enddemo
\enddemo

We are now ready to prove Theorem (2.1).

\demo{Proof of Theorem (2.1)}
Choose a system $y_1, y_2, \cdots, y_d$ of
parameters for $A$ that forms a
strong
$d$-sequence in $A$ (this choice is
possible; cf. Proposition (2.8)). Let
$x_i = y_i^2$
$(1 \leq i
\leq d)$. Then the sequence $x_1, x_2,
\cdots, x_d$ is still a strong
$d$-sequence in $A$. If
$\roman{e}(A) > 1$, then by Corollary
(2.11)
$I^2 = QI$ for the parameter ideals $Q =
(x_1^{n_1}, x_2^{n_2}, \cdots,
x_d^{n_d})$ with $n_i \geq 1$. Suppose
that $\roman{e}(A) = 1$. Then $A$ is a
regular local ring, since $A$ is unmixed,
i.e., $\dim \widehat{A}/\frak{p} = d$ for
all $\frak{p} \in \Ass \widehat{A}.$ Hence
$I^2 = QI$ by Theorem (1.1) since $Q
\subseteq
\frak{m}^2$, which  completes the proof
of Theorem (2.1). 
\qed
\enddemo

Since every parameter ideal
$\widehat{Q}$ in
$\widehat{A}$ has the form
$\widehat{Q} = Q\widehat{A}$ with $Q$ a
parameter ideal in $A$, from Theorem (2.1)
we readily get the following. 

\proclaim{Corollary (2.13)}
Let $A$ be a Noetherian local ring with
$d = \dim A \geq 2$. Assume that $A$
is unmixed, that is $\dim
\widehat{A}/\frak{p} = d$ for all
$\frak{p} \in \Ass \widehat{A}.$ Then
$A$ contains infinitely many parameter
ideals $Q$, for which the equality $I^2 =
QI$ holds true, where $I = Q : \frak{m}.$
\endproclaim

Let $A$ be a Noetherian local ring with
$d = \dim A \geq 1$. Then we say that $A$
is a generalized Cohen-Macaulay ring (or
simply,
$A$ has $FLC$), if all the local
cohomology modules
$\roman{H}_{\frak{m}}^i(A)$ $(i \ne d)$
are finitely generated $A$-modules. This
condition is equivalent to saying that
there exists an integer $\ell \gg 0$ such
that every system of parameters contained
in $\frak{m}^{\ell}$ forms a $d$-sequence
(\cite{CST}). Consequently, when $A$ is a
generalized Cohen-Macaulay ring, every
system of parameters contained in
$\frak{m}^{\ell}$ forms a strong
$d$-sequence in any order, so that by
Corollary (2.11) our local ring $A$
contains numerous parameter ideals $Q$
for which the equality $I^2 = QI$ holds
true, unless $\roman{e}(A) = 1$.
Nevertheless, even though $A$ is a
generalized Cohen-Macaulay ring with
$\roman{e}(A) > 1$, it remains subtle
whether $I^2 = QI$ for every parameter
ideal $Q$ contained in $\frak{m}^{\ell}$
$(\ell \gg 0)$. In the next section we
shall study this problem in the
case where $A$ is a Buchsbaum ring.

\head 3. Buchsbaum local rings.
\endhead
Let $A$ be a Noetherian local ring and $d
= \dim A \geq 1$. Then $A$ is said to be
a Buchsbaum ring, if the difference
$$\roman{I}(A) = \ell_A(A/Q) -
\roman{e}_Q^{0}(A)$$
is independent of the particular choice
of a parameter ideal $Q$ in $A$ and is an
invariant of $A$, where
$\roman{e}_Q^{0}(A)$ denotes the
multiplicity of $A$ with respect to $Q$.
The condition is equivalent to
saying that every system $x_1,
x_2,
\cdots, x_d$ of parameters for $A$
forms a weak $A$-sequence, that is the
equality
$$(x_1,\cdots,
x_{i-1}) : x_i = (x_1,\cdots,
x_{i-1}) : \frak{m}$$
holds true for all $1 \leq
i \leq d$ (cf. \cite{SV1}). Hence every
system of parameters for a Buchsbaum
local ring forms a strong $d$-sequence in
any order. Cohen-Macaulay local rings
$A$ are Buchsbaum rings with
$\roman{I}(A) = 0,$ and vice versa. In
this sense the notion of Buchsbaum ring
is a natural generalization of that of
Cohen-Macaulay ring.

If $A$ is a
Buchsbaum ring, then all the local
cohomology modules
$\roman{H}_{\frak{m}}^i(A)$ $(i
\ne d)$ are killed by the maximal ideal
$\frak{m}$, and one has the equality
$$\roman{I}(A) =
\sum_{i=0}^{d-1} {d-1\choose i} h^i(A)$$
where $h^i(A) =
\ell_A(\roman{H}_{\frak{m}}^i(A))$ for $0
\leq i \leq d-1$ (cf. \cite{SV2, Chap. I,
(2.6)}). It was proven by
\cite{G1, Theorem (1.1)} that for given
integers
$d$ and
$h_i \geq 0$ $(0 \leq i \leq d-1)$ there
exists a Buchsbaum local ring
$(A,\frak{m})$ such that $\dim A = d$
and
$h^i(A) = h_i$ for all $0 \leq i \leq
d-1$. One may also choose the Buchsbaum
ring $A$ so that $A$ is an integral
domain (resp. a normal domain), if $h_0 =
0$ (resp. $d \geq 2$ and $h_0 = h_1 = 0)$.
See the book \cite{SV2} for the basic
results on Buchsbaum rings and modules. 

Let $A$ be a Buchsbaum local ring with $d
=
\dim A
\geq 1$ and let
$$\roman{r}(A) = \sup_Q~\ell_A((Q :
\frak{m})/Q)$$
where $Q$ runs over parameter ideals in
$A$. Then one has the equality 
$$\roman{r}(A) =
\sum_{i=0}^{d-1}{d\choose i}h^i(A) +
\mu_{\hat{A}}(\roman{K}_{\hat{A}})$$
where $\roman{K}_{\hat{A}}$ denotes the
canonical module of $\widehat{A}$ (cf.
\cite{GSu, Theorem (2.5)}). In particular
$\roman{r}(A) <
\infty$.

We need the following, which is
implicitly known by \cite{GSu}. We
note a sketch of proof
for the sake of completeness.

\proclaim{Proposition (3.1)} 
Let $A$ be a
Buchsbaum local ring with
$d = \dim A \geq 2$. Then one
has the inequality
$\roman{r}(A/(a)) \leq \roman{r}(A)$ for every $a \in
\frak{m}$ such that $\dim A/(a) = d - 1.$
\endproclaim

\demo{Proof}
Let $B = A/(a)$. Then since
$\frak{m}{\cdot}[(0) : a] = (0)$,
from the exact sequence $$0 \to (0) : a
\to A \overset{a}\to{\to} A \to B
\to 0$$ we get a long
exact sequence
$$
\align
0 \to (0) : a 
&\to
\roman{H}_{\frak{m}}^0(A)
\overset{a}\to{\to}
\roman{H}_{\frak{m}}^0(A) \to
\roman{H}_{\frak{m}}^0(B) \\
&\to
\roman{H}_{\frak{m}}^1(A)
\overset{a}\to{\to}
\roman{H}_{\frak{m}}^1(A) \to
\roman{H}_{\frak{m}}^1(B) \\
& \cdots \\ 
&\to
\roman{H}_{\frak{m}}^i(A)
\overset{a}\to{\to}
\roman{H}_{\frak{m}}^i(A)\to
\roman{H}_{\frak{m}}^i(B) \\
& \cdots \\
&\to
\roman{H}_{\frak{m}}^d(A)
\overset{a}\to{\to}
\roman{H}_{\frak{m}}^d(A)\to 
\roman{H}_{\frak{m}}^d(B) \to
\cdots
\endalign
$$
of local cohomology modules, which
splits into the following short exact
sequences

$$0 \to \roman{H}_{\frak{m}}^i(A)\to
\roman{H}_{\frak{m}}^i(B) \to
\roman{H}_{\frak{m}}^{i+1}(A) \to 0 \ \
\  (0
\leq i \leq d-2) \ \ \ \ \text{and} \tag 3.2$$
$$0
\to
\roman{H}_{\frak{m}}^{d-1}(A)
\to
\roman{H}_{\frak{m}}^{d-1}(B)
\to [(0)
:_{\roman{H}_{\frak{m}}^{d}(A)} a] \to
0, \tag 3.3$$ because $a{\cdot}
\roman{H}_{\frak{m}}^{i}(A) =(0)$ for all
$i \ne d$. Hence $h^i(B) = h^i(A) +
h^{i+1}(A)$ $(0 \leq i \leq d-2)$ by
(3.2). Apply the functor
$\roman{Hom}_A(A/\frak{m},*)$ to 
sequence (3.3) and we have the exact
sequence
$$0 \to \roman{H}_{\frak{m}}^{d-1}(A)
\to [(0)
:_{\roman{H}_{\frak{m}}^{d-1}(B)}
\frak{m}]
\to [(0)
:_{\roman{H}_{\frak{m}}^{d}(A)}
\frak{m}].\tag 3.4$$
Hence
$$
\align
\roman{r}(B) &=
\sum_{i=0}^{d-2}{d-1\choose
i}h^{i}(B) +
\mu_{\hat{B}}(\roman{K}_{\hat{B}})\\
&=
\sum_{i=0}^{d-2}{d-1\choose
i}\{h^{i}(A) + h^{i+1}(A)\} +
\mu_{\hat{B}}(\roman{K}_{\hat{B}})  \ \
\  \\ &=
\left\{\sum_{i=0}^{d-1}{d\choose
i}h^{i}(A) - h^{d-1}(A)\right\} +
\mu_{\hat{B}}(\roman{K}_{\hat{B
}})\\
&\leq
\left\{\sum_{i=0}^{d-1}{d\choose
i}h^{i}(A) - h^{d-1}(A)\right\} +
\left\{h^{d-1}(A) + 
\mu_{\hat{A}}(\roman{K}_{\hat{A}})\right\}
 \ \ (\text{by} \ \ (3.4))\\ &=
\sum_{i=0}^{d-1}{d\choose
i}h^{i}(A)+
\mu_{\hat{A}}(\roman{K}_{\hat{A}})\\
&= \roman{r}(A)
\endalign
$$
as is claimed.
\qed
\enddemo

For the rest of this section, otherwise
specified, let $A$ be a Buchsbaum local
ring and $d = \dim A \geq 1$. Let $W =
\roman{H}_{\frak{m}}^0(A) ~(= (0) :
\frak{m})$.

To begin with
we note the following.

\proclaim{Lemma (3.5)}
Let $x_1, x_2, \cdots, x_d$ be a system of
parameters for $A$. Let $n_i \geq 1$ be
integers and put $Q = (x_1^{n_1},
x_2^{n_2}, \cdots, x_d^{n_d})$. Then 
$(Q + W) : \frak{m} = Q : \frak{m}$
if $n_i \geq 2$ for some $1 \leq i \leq
d$.
\endproclaim

\demo{Proof}
We may assume $n_d \geq 2$. Let
$L = (x_1^{n_1}, \cdots,
x_{d-1}^{n_{d-1}})$ and
$x = x_d$. Then $L : x^2 = L : x =
L :
\frak{m}$ and $xW = (0)$, since $A$ is
a Buchsbaum ring. Hence
$(Q+W):\frak{m} = Q : \frak{m}$ by Lemma
(2.4), because
$W = (0) : \frak{m} \subseteq Q
: \frak{m}$.
\qed
\enddemo

\proclaim{Theorem (3.6)}
Let $x_1, x_2, \cdots, x_d$ be a system
of parameters for $A$ and put $Q = (x_1^{n_1}, x_2^{n_2}, \cdots,
x_d^{n_d})$
with $n_i \geq 1$ $(1 \leq i \leq d)$. Let
$I = Q : \frak{m}$. Then $I^2 = QI$ if
$\roman{e}(A) > 1$ and $n_i \geq 2$ for
some $1 \leq i \leq d$.
\endproclaim

\demo{Proof}
Let $n_d \geq 2$. By Corollary
(2.11) we may assume that $d \geq 2$ and
that our assertion holds true for $d-1$.
Let $a_i = x_i^{n_i}$ $(1 \leq i \leq d)$
and put
$\bar{A} = A/(a_1)$. Then $x_2, \cdots,
x_d$ forms a system of parameters in
the Buchsbaum local ring $\bar{A}$.
Because $\roman{e}(\bar{A}) > 1$
and $n_d \geq 2$, by the
hypothesis on $d$ we get that
$\bar{I}^2 =
(\bar{a_2}, \cdots,
\bar{a_d})\bar{I}$ in $\bar{A}$, where
$\bar{a_i}$ denotes the reduction of
$a_i$ mod $(a_1)$ and $\bar{I} =
I/(a_1)$. Hence $I^2 \subseteq (a_2,
\cdots, a_d)I + (a_1)$. Since $(Q + W)
: \frak{m} = I$ by Lemma (3.5), similarly
as in the proof of Claim (2.12) we get
$(a_1) \cap I^2 = a_1I$, whence $I^2 =
QI$ as is claimed.
\qed
\enddemo

In Corollary (2.11) one needs the
assumption that $n_i \geq 2$ for $all$ $1
\leq i \leq d$. In contrast, if $A$ is a
Buchsbaum local ring, that is the case of
Theorem (3.6), this assumption is weakened
so that $n_i \geq 2$ for $some$ $1 \leq i
\leq d$. Unfortunately the assumption in
Theorem (3.6) is in general not
superfluous, as we will show in Sections
4 and 5.

The following is an immediate consequence
of Theorem (3.6).

\proclaim{Corollary (3.7)}
Let $(R,\frak{n})$ be a Buchsbaum local
ring with $\dim R \geq 2$ and
$\roman{e}(R) > 1$. Choose $f \in
\frak{n}$ so that $\dim R/(f) = \dim R -
1$ and put $A = R/(f^n)$ with $n \geq 2$.
Then the equality $I^2 = QI$ holds true
for every parameter ideal $Q$ in $A$,
where $I = Q :
\frak{m}$.
\endproclaim

Let us note one more consequence.

\proclaim{Corollary (3.8)}
Let $x_1, x_2, \cdots, x_d$ be a system
of parameters in a Buchsbaum local ring
$A$ with $d = \dim A \geq 2$ and let $Q =
(x_1^{n_1}, x_2^{n_2}, \cdots,
x_d^{n_d})$ with $n_i \geq 1$ $(1 \leq i
\leq d)$. Then $I^2 = QI$ if $n_i, n_j
\geq 2$ for some $1 \leq i, j \leq d$ with
$i \ne j$.
\endproclaim

\demo{Proof}
Thanks to Theorem (3.6) we may assume
that $\roman{e}(A) = 1$. Let $B = A/W$. Then $B$ is
a regular local ring with $\dim B = d \geq
2$, because $\roman{e}(B) = 1$ and $B$
is unmixed (cf.
\cite{CST}). We have $\ell_B((QB +
\frak{m}^2B)/\frak{m}^2B) \leq d - 2$, since $x_i^{n_i},
x_j^{n_j}
\in
\frak{m}^2$. Therefore $(IB)^2 =
(QB){\cdot}(IB)$ by Theorem (1.1), because $IB = QB :
\frak{m}B$ (recall that $I = (Q + W) :
\frak{m}$ by Lemma (3.5)). Hence
$I^2
\subseteq QI + W$, so that we have  $I^2 =
QI$ since $W \cap Q = (0)$ (cf. (2.6) and (2.7)).
\qed
\enddemo

We now turn to other topics.

\proclaim{Theorem (3.9)}
Let $A$ be a Buchsbaum local ring with $d
=
\dim A \geq 1$ and $\roman{e}(A) > 1$.
Let $Q$ be a parameter ideal in $A$ and
put $I = Q : \frak{m}$. Then $I^2 = QI$
if $\ell_A(I/Q) = 
\roman{r}(A)$.
\endproclaim

\demo{Proof}
Let $W = \roman{H}_{\frak{m}}^0(A)$. Then
$\frak{m}W = (0)$ and $Q \subseteq Q + W
\subseteq I \subseteq (Q + W) :
\frak{m}$. Hence
$$\ell_A(I/Q) = \ell_A(I/(Q+W)) +
\ell_A(W)$$
because $W \cap Q = (0)$. Assume that $d = 1$. Then
$\roman{r}(A) =
\ell_A(W) +
\mu_{\hat{A}}(\roman{K}_{\hat{A}}) = \ell_A(I/Q)$. Since $A/W$
is a Cohen-Macaulay ring and $\roman{H}_{\frak{m}}^1(A)
\cong \roman{H}_{\frak{m}}^1(A/W)$, we
have
$$\mu_{\hat{A}}(\roman{K}_{\hat{A}}) =
\roman{r}(A/W) = \ell_A\left([(Q+W):\frak{m}]/(Q+W)\right)$$
so that 
$$\ell_A\left([(Q+W):\frak{m}]/(Q+W)\right)
= \mu_{\hat{A}}(\roman{K}_{\hat{A}}) = \ell_A(I/Q) - \ell_A(W)
= \ell_A(I/(Q+W)).$$
Hence $(Q+W) : \frak{m} = I$ and so $I^2
=QI$ (cf. Proof of Corollary (2.11)). 

Assume now that $d \geq 2$ and that our
assertion holds true for $d-1$. Let
$Q = (a_1, a_2, \cdots, a_d)$ and put
$\bar{A} = A/(a_1)$, $\bar{Q} = Q/(a_1)$,
$\bar{I} = I/(a_1)$, and $\bar{\frak{m}}
= \frak{m}/(a_1)$. Then $\bar{I} =
\bar{Q} : \bar{\frak{m}}$ and
$\roman{r}(\bar{A}) \geq
\ell_{\bar{A}}(\bar{I}/\bar{Q}) =
\ell_A(I/Q) = \roman{r}(A)$. Hence by
Proposition (3.1) we get
$\roman{r}(\bar{A}) =
\ell_{\bar{A}}(\bar{I}/\bar{Q})$, so that
$\bar{I}^2 = \bar{Q}~\bar{I}$ by the
hypothesis on $d$. Thus $I^2 \subseteq
(a_2,, \cdots, a_d)I + (a_1)$ and then
the equality $I^2 = QI$ follows similarly
as in the proof of Claim (2.12).
\qed
\enddemo

The following is a direct consequence of
Theorem (3.9), which may account well for
the reason why $I^2 = QI$ in
Cohen-Macaulay rings $A$.

\proclaim{Corollary (3.10)}Let $A$ be a
Buchsbaum local ring with $d = \dim A
\geq 1$ and assume that the index 
$\ell_A((Q:\frak{m})/Q)$ of
reducibility of $Q$ is independent of
the choice of a parameter ideal
$Q$ in $A$.  If
$\roman{e}(A) > 1$,  then the equality
$I^2 = QI$ holds true for every parameter
ideal $Q$ in $A$, where $I = Q :
\frak{m}$.
\endproclaim

The hypothesis of Corollary (3.10) may
be satisfied even though $A$ is not a
Cohen-Macaulay ring. Let $B = \Bbb{C}[[X,
Y, Z]]/(Z^2 - XY)$ where $\Bbb{C}[[X,
Y, Z]]$ is the formal power series ring
over the field $\Bbb{C}$ of complex numbers, and
put $$A = \Bbb{R}[[x, y, z, ix, iy,
iz]]$$
where $\Bbb{R}$ is the field of real
numbers, $i = \sqrt{-1}$, and $x, y$, and
$z$ denote the reduction of $X,Y$, and
$Z$ mod $(Z^2-XY)$. Then $A$ is a
Buchsbaum local integral domain with
$\dim A= 2$, $\depth A = 1$, and
$\roman{e}(A) = 4$. For this ring $A$ one
has the equality $$\ell_A((Q :
\frak{m})/Q) = 4$$ for every parameter
ideal $Q$ in
$A$ (\cite{GSu, Example (4.8)}). Hence
by Corollary (3.10),
$I^2 = QI$ for all parameter ideals $Q$ in
$A$. 

The following theorem (3.11) gives an
answer to the question raised in the
previous section.  The authors know
no example of Buchsbaum local rings
$A$ with
$\roman{e}(A) > 1$ such that $I^2 \ne
QI$ for some parameter ideal $Q \subseteq
\frak{m}^2$.

\proclaim{Theorem (3.11)}
Let $A$ be a Buchsbaum local ring and
assume that $\dim A \geq 2$ or that
$\dim A = 1$ and $\roman{e}(A) >1$. Then
there exists an integer $\ell \gg 0$
such that $I^2 = QI$ for every parameter
ideal
$Q \subseteq \frak{m}^{\ell}$.
\endproclaim

To 
prove this theorem we need one more
lemma. Let $A$ be an arbitrary Noetherian
local ring with the maximal ideal
$\frak{m}$ and $d = \dim A \geq 1$. Let
$f : M \to N$ be a homomorphism of
$A$-modules. Then we say that $f$ is
surjective (resp. bijective) on the socles, if the induced
homomorphism
$$f_{*} : \roman{Hom}_A(A/\frak{m},M) =
(0) :_M \frak{m}
\to
\roman{Hom}_A(A/\frak{m},N)= (0)
:_N \frak{m}$$ is an
epimorphism (resp. an isomorphism). 

Let $Q = (a_1, a_2, \cdots, a_d)$ be a
parameter ideal in $A$ and let $M$ be an
$A$-module. For each integer $n \geq 1$
we denote by $\underline{a}^n$ the
sequence $a_1^n, a_2^n, \cdots,
a_d^n$. Let
$\roman{K}_{\bullet}(\underline{a}^n)$ be the
Koszul complex of $A$ generated by the
sequence $\underline{a}^n$ and let
$$\roman{H}^{\bullet}(\underline{a}^n;M) =
H^{\bullet}(\roman{Hom}_A(\roman{K}_{\bullet}
(\underline{a}^n),M))$$
be the Koszul cohomology module of
$M$. Then for every $p \in \z$ the family
$\{\roman{H}^p(\underline{a}^n;M)\}_{n
\geq 1}$ naturally forms an inductive
system of
$A$-modules, whose limit
$$\roman{H}_{\underline{a}}^p(M) =
\underset{n
\to
\infty}\to{\lim}\roman{H}^p
(\underline{a}^n;M)$$
is isomorphic to the local cohomology
module
$$\roman{H}_{\frak{m}}^p(M) =\underset{n
\to
\infty}\to{\lim}
\roman{Ext}_A^p(A/\frak{m}^n,M).$$
For each $n \geq 1$  and $p \in \z$ let 
$\varphi_{\underline{a},M}^{p,n} :
\roman{H}^p(\underline{a}^n;M)\to
\underset{n
\to
\infty}\to{\lim}\roman{H}^p_{\underline{a}}(M)$ denote the
canonical homomorphism into the limit. With this notation we
have the following.

\proclaim{Lemma (3.12)}
Let $A$ be a Noetherian local ring with
the maximal ideal $\frak{m}$ and $d =
\dim A \geq 1$. Let $M$ be a finitely
generated $A$-module. Then there exists
an integer $\ell \gg 0$ such that for
all systems $a_1, a_2, \cdots, a_d$ of
parameters for $A$ contained in
$\frak{m}^{\ell}$ and for all 
$p \in \z$ the canonical homomorphisms
$$\varphi_{\underline{a},M}^{p,1} :
\roman{H}^p(\underline{a};M) \to
\roman{H}_{\underline{a}}^p(M) =
\underset{n \to
\infty}\to{\lim}\roman{H}^p
(\underline{a}^n;M)$$ into the inductive limit are surjective
on the socles.
\endproclaim

\demo{Proof}
First of all, choose $\ell \gg 0$ so that
the canonical homomorphisms
$$\varphi_{\frak{m},M}^{p,\ell} :
\roman{Ext}_A^p(A/\frak{m}^\ell, M) \to
\roman{H}_{\frak{m}}^p(M) =
\underset{n \to
\infty}\to{\lim}\roman{Ext}_A^p
(A/\frak{m}^n, M)$$ are surjective on
the socles for all $p \in \z$. This
choice is possible, because
$\roman{H}_{\frak{m}}^p(M) = (0)$ for
almost all $p \in \z$ and the socle of
$[(0):_{\roman{H}_{\frak{m}}^p(M)}
\frak{m}]$ of
$\roman{H}_{\frak{m}}^p(M)$ is finitely
generated. Let $Q=(a_1, a_2, \cdots,
a_d)$ be a parameter ideal in $A$ and
assume that $Q \subseteq
\frak{m}^{\ell}$. Then, since $\sqrt{Q}
= \sqrt{\frak{m}^{\ell}} = \frak{m}$,
there exists an isomorphism
$\theta_M^p : \roman{H}_{\frak{m}}^p(M)
\to
\roman{H}_{Q}^p(M) = \underset{n \to
\infty}\to{\lim}\roman{Ext}_A^p
(A/Q^n, M)$ which makes
the diagram
$$\CD
\roman{Ext}_A^p
(A/\frak{m}^\ell, M)
@>\varphi_{\frak{m},M}^{p, \ell}>>
\roman{H}_{\frak{m}}^p (M)\\
\alpha@VVV @VV{\theta_M^p}V \\
\roman{Ext}_A^p
(A/Q, M) @>>\varphi_{Q,M}^{p,1}>
\roman{H}_{Q}^p (M)
\endCD
$$
commutative, where the vertical map
$\alpha : \roman{Ext}_A^p
(A/\frak{m}^{\ell}, M) \to \roman{Ext}_A^p
(A/Q, M)$ is the homomorphism induced
from the epimorphism $A/Q \to
A/\frak{m}^{\ell}$. Hence the
homomorphism $\varphi_{Q,M}^{p,1}$ is
surjective on the socles, since so is
$\varphi_{\frak{m},M}^{p,\ell}$.  Let $n
\geq 1$ be an integer and let
$$\cdots \to F_i\to \cdots \to
F_1 \to F_0 = A \to A/Q^n \to 0$$ be a
minimal free resolution of $A/Q^n$. Then
since $(\underline{a}^n) \subseteq Q^n$,
the epimorphism $$\varepsilon :
A/(\underline{a}^n) \to A/Q^n$$ can be
lifted to a homomorphism of complexes:
$$\CD
\cdots @>>> F_i @>>> \cdots @>>>
F_1 @>>> F_0 = A @>>> A/Q^n @>>> 0 \\
@. @AAA @. @AAA @| @AA{\varepsilon}A \\
\cdots @>>> K_i @>>> \cdots @>>>
K_1 @>>> K_0 = A @>>> A/(\underline{a}^n)
@>>> 0
\endCD
$$
where $K_{\bullet} =
\roman{K}_{\bullet}(\underline{a}^n)$.
Taking the $M$-dual of these two
complexes and passing to the cohomology
modules, we get the natural homomorphism
$$\alpha_{M}^{p,n} :
\roman{Ext}_A^p(A/Q^n,M) \to
\roman{H}(\underline{a}^n;M)$$
$(p \in \z, n\geq 1)$ of inductive
systems, whose limit $$\alpha^p_M : \roman{H}_Q^p(M) \to
\roman{H}_{\underline{a}}^p(M)$$ is necessarily an
isomorphism
for all $p \in \z$. Consequently, thanks
to the commutative diagram
$$\CD
\roman{Ext}_A^p
(A/Q, M)
@>\varphi_{Q,M}^{p,1}>>
\roman{H}_{Q}^p (M)\\
\alpha_{M}^{p,1}@VVV @VV{\alpha^p_M}V \\
\roman{H}^p(\underline{a};M)
@>>\varphi_{\underline{a},M}^{p,1}>
\roman{H}_{\underline{a}}^p (M)
\endCD
$$
we get that for all $p \in \z$ the
homomorphism 
$$\varphi_{\underline{a},M}^{p,1}
: \roman{H}^p(\underline{a};M)
\to \roman{H}_{\underline{a}}^p (M)$$
is surjective on the socles, because so is
$\varphi_{Q,M}^{p,1}$.
\qed
\enddemo

\proclaim{Corollary (3.13)}
Let $A$ be a Buchsbaum local ring with $d = \dim A \geq 1$.
Then there exists an integer $\ell \gg 0$ such that the index
$\ell_A((Q : \frak{m})/Q)$ of reducibility of $Q$ is
independent of $Q$ and equals $\roman{r}(A)$ for all
parameter ideals $Q  \subseteq \frak{m}^\ell$. 
\endproclaim

\demo{Proof}
Choose an integer $\ell \gg 0$ so that
the canonical homomorphism
$$\varphi_{\underline{a},A}^{d,1}
: A/Q = \roman{H}^d(\underline{a};A)
\to \roman{H}_{\underline{a}}^d (A)$$ is surjective on the socles for
every parameter ideal $Q = (a_1, a_2, \cdots,
a_d) \subseteq
\frak{m}^{\ell}$. Then since $A$ is a
Buchsbaum local ring, we get that
$$\roman{Ker}~
\varphi_{\underline{a},A}^{d,1} =
\sum_{i=1}^d\left[\left((a_1, \cdots,
\overset{\vee}\to{a_i}, \cdots, a_d) :
a_i\right) + Q\right]/Q$$
(\cite{G2,
Theorem
(4.7)}), $\frak{m}{\cdot}[\roman{Ker}~
\varphi_{\underline{a},A}^{d,1}]=(0)$, and
$\ell_A(\roman{Ker}~
\varphi_{\underline{a},A}^{d,1}) =
\sum_{i=0}^{d-1}{d\choose i}h^i(A)$
(\cite{G2, Proposition (3.6)}). Because
$\mu_{\hat{A}}(\roman{K}_{\hat{A}}) =
\ell_A((0)
:_{\roman{H}^d_{\underline{a}}(A)}
\frak{m})$, the surjectivity of the homomorphism
$\varphi_{\underline{a},A}^{d,1}$ on the socles guarantees
that
$$\ell_A(I/Q)=\sum_{i=0}^{d-1}{d\choose
i}h^i(A)+\mu_{\hat{A}}(\roman{K}_{\hat{A}})$$
where
$I = Q :
\frak{m}$. Hence
$\roman{r}(A) = \ell_A(I/Q)$.
\qed
\enddemo

We are now ready to prove Theorem (3.11). 

\demo{Proof of Theorem (3.11)}
Thanks to Theorem (3.9) and Corollary
(3.13) we may assume that
$\roman{e}(A) = 1$ and $d
\geq 2$. Let $W =
\roman{H}^0_{\frak{m}}(A)$ and $B =
A/W$. Then $B$ is a regular local ring
with $d = \dim B \geq 2$. We choose a
parameter ideal $Q$ in $A$ so that $Q
\subseteq
\frak{m}^2$. Let $J = QB : \frak{m}B$.
Then since $QB \subseteq (\frak{m}B)^2$,
by Theorem (1.1) we get $J^2 =
QB{\cdot}J$. Because $B/QB$ is a
Gorenstein ring and $QB \subseteq IB
\subseteq J$, we have either $IB = QB$ or
$IB =J$. In any case $I^2 \subseteq QI +
W$, so that $I^2 =QI$, because $W \cap Q =
(0)$. 
\qed
\enddemo

\head 4. Evaluation of $\roman{r}_Q(I)$ in the case where
$\dim A = 1$. \endhead

In this section let $A$ be a Buchsbaum local ring and assume
that $\dim A = 1$. Let $W = \roman{H}_{\frak{m}}^0(A)~(=
(0) :
\frak{m})$ and
$ e =
\roman{e}(A)$. Then $\roman{r}(A) = \ell_A(W) +
\roman{r}(A/W)$ and $\roman{r}(A/W) \leq \max \{1, e-1\}$,
since $A/W$ is a Cohen-Macaulay local ring with $
\roman{e}(A/W) = e$ (cf.
\cite{HK, Bemerkung 1.21 b)}). The
purpose is to prove the following. 

\proclaim{Theorem (4.1)}
Suppose that~$e > 1$. Let $Q$ be a
parameter ideal in
$A$ and put $I = Q : \frak{m}$. Then
$$\roman{r}_Q(I) \leq \roman{r}(A) - \ell_A(W) + 1 =
\roman{r}(A/W) - \ell_A(I/(Q + W)) + 1.$$

\endproclaim

\demo{Proof}
Let $Q = (a)$ and put $I_n = I^{n+1} : a^n$ ($n
\geq 0$). Then $I_0 = I$ and $I_n \subseteq I_{n+1}$. We
have $I_n \subseteq (Q+W):\frak{m}$. In fact, let $x \in
I_n$ and $\alpha \in \frak{m}$. Then
$a^n(\alpha x) \in \frak{m}I^{n+1} \subseteq (a^{n+1})$ by
Proposition (2.3). Let $a^n(\alpha x) = a^{n+1}y$ with $y
\in A$. Then $\alpha x - ay \in (0) : a^n = W$, whence $x
\in (Q+W) : \frak{m}$. We furthermore have the following.

\proclaim{Claim (4.2)}
Let $n \geq 0$ and assume that $I_n = I_{n+1}$. Then
$I^{n+2} = QI^{n+1}$.
\endproclaim

\demo{Proof of Claim (4.2)}
Let $x \in I^{n+2} \subseteq (a^{n+1})$ and write $x =
a^{n+1}y$ with $y \in A$. Then $y \in I^{n+2} : a^{n+1} =
I_n$, so that $x = a(a^ny) \in QI^{n+1}$. Thus $I^{n+2} =
QI^{n+1}$.
\qed
\enddemo

Let $\ell = \ell_A(I/(Q+W))$. Then
$\roman{r}(A/W) = \ell_A([(Q + W) :
\frak{m}]/(Q+W))
\geq 
\ell$. Since
$\ell_A(I/Q) =
\ell_A(I/(Q+W)) + \ell_A(W)$ (cf. Proof of Theorem (3.9)),
we get
$$
\align
\roman{r}(A) - \ell_A(I/Q) + 1 &= [\roman{r}(A/W) +
\ell_A(W)] - [\ell_A(I/(Q+W)) + \ell_A(W)] + 1 \\
&=
\roman{r}(A/W) - \ell_A(I/(Q+W)) + 1 \\
&= \roman{r}(A/W) - \ell + 1.
\endalign
$$
Assume that $\roman{r}_Q(I) > \roman{r}(A/W) - \ell + 1$
and put $n = \roman{r}(A/W) - \ell + 2$. Then
$\roman{r}_Q(I) \geq n \geq 2$, so that by Claim (4.2) $I_i
\ne I_{i+1}$ for all $0 \leq i \leq n - 2$. Hence we have
a chain
$$Q + W \subseteq I_0 = I \subsetneq I_1 \subsetneq \cdots
\subsetneq I_{n - 2} \subsetneq I_{n - 1} \subseteq (Q+W) :
\frak{m}$$
of ideals, so that $\roman{r}(A/W)
=
\ell_A([(Q+W):\frak{m}]/(Q+W)) \geq (n - 1 ) + \ell =
\roman{r}(A/W) + 1$, which is absurd. Thus $\roman{r}_Q(I)
\leq \roman{r}(A/W) - \ell + 1$.
\qed
\enddemo

Suppose that $e > 1$ and let $Q$ be a parameter ideal in
$A$. Let $I = Q : \frak{m}$. Then $I \supseteq Q + W$. We
have by Theorem (4.1) that
$\roman{r}_Q(I)
\leq
\roman{r}(A/W) \leq e - 1$, if $I \supsetneq Q+W$. If $I =
Q + W$, then $I^2 = Q^2$ because $\frak{m}W = (0)$, so
that $I^n = Q^n$ for all $n \geq 2$. Thus we have

\proclaim{Corollary (4.3)}
Let $A$ be a Buchsbaum local ring with $\dim A = 1$ and $e
= \roman{e}(A) > 1$. Then 
$$\sup_Q~\roman{r}_Q(Q:\frak{m}) \leq
e-1$$
where $Q$ runs over parameter ideals in
$A$.
\endproclaim

The evaluations in
Theorem (4.1) and Corollary (4.3) are
sharp, as we shall show in the following
example. The example shows that for every
integer $e \geq 3$ there exists a
Buchsbaum local ring $A$ with $\dim A =
1$ and $\roman{e}(A) = e$ which contains
a parameter ideal
$Q$ such that $\roman{r}_Q(I) = e-1$,
where $I = Q : \frak{m}$. Hence the
equality
$I^2 =QI$ fails in general to hold, even
though $A$ is a Buchsbaum local ring with
$\roman{e}(A) >1$. The reader may consult
the forthcoming paper
\cite{GSa} for higher-dimensional
examples of higher depth.

Let $k$ be a field and $3 \leq e \in \z$. Let $S=k[X_1,
X_2, \cdots, X_e]$ and
$P=k[t]$ be the polynomial rings over $k$. We regard $S$ and
$P$ as
$\z$-graded rings whose gradings are given by $S_0 =k$,
$S_{e + i - 1} \ni X_i$ $(1 \leq i \leq e)$ and $P_0=k$,
$P_1 \ni t$. Hence $S_n = (0)$ for
$1 \leq n \leq e$, where $S_n$
denotes the homogeneous component of
$S$ with degree $n$. Let
$\varphi : S
\to P$ be the $k$-algebra map defined
by
$\varphi(X_i) = t^{e + i - 1}$ for all
$1 \leq i \leq e$. Then $\varphi$ is a
homomorphism of graded rings, whose image is the semigroup
ring $k[t^e, t^{e+1}, \cdots, t^{2e-1}]$, and whose kernel
$ \frak{p}$ is minimally generated by the 2 by 2 minors of
the matrix
$$\Bbb{M} = \left(\matrix X_1 & X_2 & \cdots & X_{e-1} & X_e
\\
X_2 & X_3 & \cdots & X_{e} & X_1^2 
\endmatrix \right).
$$
Let $\Delta_{ij}$ $(1 \leq i, j \leq e)$ be the determinant
of the matrix consisting of the $i~^{\underline{th}}$ and 
$j~^{\underline{th}}$ columns of $\Bbb{M}$, that is
$$\Delta_{ij} = \left| \matrix X_i &
X_j \\ X_{i+1} & X_{j+1}
\endmatrix
\right|,
$$ where $X_{e+1}= X_1^2$ for convention. We put $\Delta =
\Delta_{2,e}$ and let $N = 
S_+ ~(= \bigoplus_{n\geq 1}S_n)$, the
unique graded maximal ideal in $S$. Let
$$\frak{a} = \left(\Delta_{ij} \mid 
1 \leq i < j \leq e\ \ \text{such
that} \ \   (i,j)
\ne (2,e)\right) + \Delta N$$ and put $R = S/\frak{a}$, $M
= R_+$, $A= R_M$, and $\frak{m} = MA$.
Let
$x_i = X_i$ mod
$\frak{a}$ $(1 \leq i \leq e)$ and $\delta = \Delta$ mod
$\frak{a}$. We then have the following.

\proclaim{Lemma (4.4)}
$\dim R = 1$, $\roman{H}_M^0(R) = (\delta) \ne (0)$, and
$M\delta = (0)$.
\endproclaim

\demo{Proof}
We certainly have $M\delta = (0)$. Look at the
canonical exact sequence
$$0 \to \frak{p}/\frak{a} = (\delta) 
\to R
\to S/\frak{p} \to 0, \tag 4.5$$
where $\frak{p} = \Ker
\varphi$.
Then, since $M\delta = (0)$ and $S/\frak{p} =k[t^e,
t^{e+1}, \cdots, t^{2e-1}]$ is a 
Cohen-Macaulay integral domain with
$\dim S/\frak{p} = 1$, we get that 
$\dim R = 1$ and $\roman{H}_M^0(R) = (\delta)$. The
assertion $\delta \ne 0$ follows from 
the fact that
$\{\Delta_{ij}\}_{1 \leq i < j \leq e}$
is a $minimal$ system
of generators for the ideal 
$\frak{p}$. 
\qed
\enddemo

Let $T = k[t^e, t^{e+1}, \cdots, t^{2e-1}]$ and $\frak{n} =
T_+$. Then $\frak{n}=(t^e, t^{e+1},
\cdots, t^{2e-1})T$ and
$\frak{n}^2 = t^e\frak{n}$. Hence
$$\roman{r}(T_{\frak{n}}) = \ell_T((t^eT : \frak{n})/t^eT) =
\ell_T(\frak{n}/t^eT) = e - 1.$$ We have $M^2 = x_1M +
(\delta)$, because
$\frak{n}^2 = t^e\frak{n}$ and $\delta \in M^2$.
Hence
$M^3 = x_1M^2$, so that $\roman{e}(A) =
\roman{e}_{x_1A}^0(A) =
\roman{e}_{x_1A}^0(T_{\frak{n}}) =
\ell_{T}(T/t^eT) = e$  (cf. (4.5)). Thus $A$ is a
Buchsbaum ring with $\dim A = 1$ and
$\roman{e}(A) = \roman{r}(A) = e$. In particular, $\delta
\not\in (x_1)$, since $(x_1) \cap
\roman{H}_M^0(R) = (0)$ (recall that
$x_1$ is a parameter of $R$).

We put $J = (x_1) : M$. 

\proclaim{Proposition (4.6)}
The following assertions hold true.
\roster
\item $J = (x_1, x_2, \delta)$.
\item $J^n = (x_1,x_2)^n$ for all $n \geq 2$.
\item $\ell_R(J/(x_1)) = 2$.
\endroster  
\endproclaim

\demo{Proof}
We firstly  notice that 
$$
\align
\frak{a} + X_1 &\supseteq (X_1) + (X_2, X_3X_e)(X_2,
\cdots, X_e) \tag 4.7 \\
&+ (\Delta_{ij} \mid 3 \leq i,j \leq e, ~~i+j = e + 2) \\
&+ (X_iX_j \mid 3 \leq i,j \leq e,~~i+j \ne e + 3).
\endalign 
$$
In fact, $\Delta \equiv -X_3X_e$ mod $(X_1)$ and
$\Delta_{1,j}=X_1X_{j+1} -X_2X_j \equiv -X_2X_j$ mod
$(X_1)$, we get $\frak{a} + (X_1) \supseteq (X_1) + (X_2, X_3X_e)(X_2,
\cdots, X_e)$. Let $3 \leq i, j \leq e$. If $i+j = e+2$,
then $(i,j) \ne (2,e)$ and $(j,i) \ne (2,e)$, so that
$\Delta_{ij} \in \frak{a}$. Assume that $i+j \ne e+3$. We
will show $X_iX_j \in \frak{a} + (X_1)$ by induction on
$i$. If $i = 3$, then $3 \leq j < e$ and $\Delta_{2j} =
X_2X_{j+1} - X_3X_j \in \frak{a}$, whence $X_3X_j \in
\frak{a} + (X_1)$, because $X_2X_{j+1} \in \frak{a} +
(X_1).$ Assume that $i \geq 4$ and that our assertion holds
true for $i - 1$. Then $3 \leq i - 1 < e$, so that
$\Delta_{i-1, j} = X_{i-1}X_{j+1} - X_iX_j \in \frak{a}$.
Hence
$X_iX_j \in \frak{a} + (X_1)$, because $X_{i-1}X_{j+1} \in
\frak{a} + (X_1)$ by the hypothesis on $i$.

Let $B = S/(\frak{a} + (X_1))$ and
$\frak{q} = B_+$. Then
$(B,\frak{q})$ is an Artinian graded local ring. For
the moment, let us denote by $y_i$ the
reduction of
$X_i$ mod
$\frak{a} + (X_1)$ $(2 \leq i \leq e)$
and by $\rho$ the reduction of
$-\Delta$ mod
$\frak{a} + (X_1)$. Hence $\frak{q} = 
(y_2, \cdots, y_e)$
and $\rho = y_3y_e$. We will check that $\frak{q}^2 =
(\rho)$. To see this, let $2
\leq i, j \leq e$ and assume that $y_iy_j \ne
0$. Then $3 \leq i,j \leq e$ and 
$i+j = e+3$ by (4.7), whence 
$y_iy_j = \rho$, because $\rho =y_3 y_e$ and
$y_\alpha y_{\beta +1} = y_{\alpha +1}y_\beta$ whenever
$3
\leq \alpha,\beta \leq e$ with $\alpha +\beta =e + 3$. Hence
$\frak{q}^2 = (\rho)$, so that $\frak{q}^3 = (0)$ because
$N{\cdot}\Delta \subseteq \frak{a}$. We have
$\rho \ne 0$, since $\Delta \not\in \frak{a} +
(X_1)$ (recall that $\delta \not\in
(x_1)$). Now let
$\varphi
\in (0) :
\frak{q}$ and write
$\varphi = c + \sum_{i=2}^e c_iy_i + d\rho$ with $c, c_i, d
\in k$. Then because $(0) : \frak{q}$ is a graded
ideal in $B$ and $c_iy_i \in B_{e+i-1}$ for $2 \leq i \leq
e$ and $\rho \in B_{3e+1}$, we get $c, c_iy_i, d\rho \in
(0) : \frak{q}$. Hence $c = 0$, because $(0) : \frak{q}
\subseteq \frak{q}$. We have
$c_i = 0$ for all $3 \leq i \leq e$, because $\rho =
y_\alpha y_{e - \alpha + 3} \ne 0$ for
all
$3 \leq \alpha \leq e$. Thus $\varphi =
c_2y_2 + d\rho \in (y_2, \rho)$. Hence
$(0) :
\frak{q} = (y_2, \rho)$ by (4.7), so
that we have $J = (x_1, x_2,
\delta)$ in
$R$. Assrtions (2) and (3) are now clear.
\qed
\enddemo

\proclaim{Theorem (4.8)}
$J^e = x_1J^{e-1}$ but $J^{e-1} \ne x_1J^{e-2}$.
\endproclaim

\demo{Proof}
Assume that $J^{e-1} = x_1J^{e-2}$. Then $J^{e-1} \ni
x_2^{e-1} = x_2^2x_2^{e-3} = x_1{\cdot}x_2^{e-3}x_3$. Let
$x_1{\cdot}x_2^{e-3}x_3 = x_1\eta$ with $\eta
\in J^{e-2}$. Then $x_2^{e-3}x_3 - \eta \in (0) : x_1 =
(\delta)$. We write 
$$x_2^{e-3}x_3 =
\eta + \delta \xi$$
with
$\xi \in R$. If $e = 3$, then $x_3 \in J = (x_1, x_2,
\delta) \subseteq (x_1, x_3^2)$, which is impossible. Hence $e \geq
4$ and so $\eta \in (x_1)$, since $\eta \in J^{e-2}
\subseteq J^2$ and $J^2 = (x_1, x_2)^2 = (x_1^2, x_1x_2,
x_2^2)
\subseteq (x_1)$ (cf. Proposition
(4.2) (2); recall that
$x_2^2 = x_1x_3$). Hence
$\delta \xi \in (x_1) \cap \roman{H}_M^0(R) = (0)$, because
$x_2^{e-3}x_3 = x_2x_3{\cdot}x_2^{e-4} = x_1x_4x_2^{e-4}
\in (x_1)$. Thus by Proposition
(4.2) (2) $$x_2^{e-3}x_3 =
\eta \in (x_1, x_2)^{e-2} =
(x_1^ix_2^{e - 2 - i}
\mid 0
\leq i
\leq e - 2). \tag 4.9$$ Here we notice
that $R= \bigoplus_{n \geq 0}R_n$ is a
graded ring and that
$\roman{deg}~(x_1^ix_2^{e - 2 - i}) =
e^2 - e - i -2$,
$\roman{deg}~(x_2^{e-3}x_3) = e^2 - e
- 1$. Then, since 
$1 \leq i + 1 = (e^2 - e - 1) - (e^2 - e - i
-2)\leq e-1$ for 
$0
\leq i
\leq e - 2$ and $R_n = (0)$ for $1 \leq  n \leq e - 1$,
by (4.9) we get
$x_2^{e-3}x_3 = 0$, whence
$X_2^{e-3}X_3 \in \frak{p} = \Ker \varphi$, which is
impossible. Thus
$J^{e-1} \ne x_1 J^{e-2}$. Since $J^e = x_1J^{e-1} +
(x_2^e)$, the equality $J^e = x_1J^{e-1}$ follows
from Corollary (4.3), or more directly
from the following.

\proclaim{Claim (4.10)}
$x_2^e = x_1^{e+1}$.
\endproclaim

\demo{Proof of Claim (4.10)}
It suffices to show $x_2^e = x_1^nx_2^{e - n - 1}x_{n+2}$
for all $1 \leq n \leq e - 2$. Since $x_2^e =
x_1x_3{\cdot}x_2^{e-2}$, the assertion is obviously true
for $n = 1$. Let $n \geq 2$ and assume that the equality
holds true for $n - 1$. Then 
$$
\align
x_2^e &=x_1^{n-1}x_2^{e-n}x_{n+1} \\
&= x_1^{n-1}x_2^{e - n -1}{\cdot}x_2x_{n+1} \\
& = x_1^nx_2^{e - n - 1}x_{n+2},
\endalign
$$
because $x_2x_{n+1} = x_1x_{n+2}$.
Hence $x_2^e = x_1^{e-2}{\cdot}x_2x_e =
x_1^{e-2}x_1^3 = x_1^{e+1}$. 
\qed
\enddemo
\enddemo

Let $Q = x_1A$ and $I = Q : \frak{m} ~(=JA)$. Then in our
Buchsbaum local ring $A$ we have $I^e = x_1I^{e-1}$ but
$I^{e-1} \ne x_1 I^{e-2}$. Because $\roman{e}(A) =
\roman{r}(A) = e$, this example shows 
the evaluations in
Theorem (4.1) and Corollary (4.3) are
really sharp.

\head 5. Examples
\endhead

In this section we shall explore two
examples. One is to show that the
equality $I^2 =QI$ may hold true for
$all$ parameter ideals $Q$ in $A$, even
though $A$ is not a generalized
Cohen-Macaulay ring. As is shown in
the previous section, the equality
$I^2 =QI$ fails in general to hold, even
though $A$ is a Buchsbaum local ring with
$\roman{e}(A) > 1$. In this section we
will also explore one counterexample of
dimension 1 and give complete criteria
of the equality 
$I^2 = QI$ for parameter ideals $Q$ in
the example.

Throughout this section let
$(R,\frak{n})$ be a 3-dimensional regular
local ring and let $\frak{n} = (X, Y,
Z)$. Firstly, let $\ell
\geq 1$ be an integer and put $$A =
R/(X^{\ell}) \cap (Y,Z).$$ Let $x, y$, and
$z$ denote the reduction of $X, Y$, and
$Z$ mod
$(X^\ell) \cap (Y,Z)=(X^\ell Y, X^\ell
Z)$. Let
$\frak{p} = (y,z)$. Then $\frak{m} = (x)
+
\frak{p}$ and $(x^\ell) \cap \frak{p} =
(0)$ in
$A$, where $\frak{m}$ denotes the
maximal ideal in $A$. Let $B =
A/(x^\ell)$.
Then there exists exact sequences 
$$0 \to A/\frak{p}
\overset{\alpha}\to{\to} A \to B
\to 0  \ \ \ \text{and} \tag 5.1$$
$$0 \to A/(x)
\overset{\beta}\to{\to} B
\to A/(x^{\ell -1})
\to 0 \tag 5.2$$
of $A$-modules, where the homomorphisms
$\alpha$ and $\beta$ are defined by
$\alpha(1) = x^{\ell}$ and $\beta(1) =
x^{\ell - 1}$ mod $(x^{\ell})$. Since
$A/\frak{p}$ is a DVR and $B$ is a
hypersurface with $\dim B = 2$, we get by
(5.1) that
$$\dim A = 2, ~~\depth A = 1,
~~\text{and}~~
\roman{H}^1_{\frak{m}}(A/\frak{p}) \cong
\roman{H}^1_{\frak{m}}(A).$$ Hence $A$
is not a generalized Cohen-Macaulay ring.
Let
$\frak{q} = (x-y,z)$. Then
$\frak{m}^{\ell + 1} =
\frak{q}\frak{m}^{\ell}$, since
$\frak{m} = (x) + \frak{q}$ and
$x^{\ell + 1} = (x-y)x^{\ell}$.
Consequently  by (5.1) we get
$$\roman{e}(A) =
\roman{e}_{\frak{q}}^0(A) =
\roman{e}_{\frak{q}}^0(B) =
\ell_A(B/\frak{q}B) =
\ell_R(R/(X^{\ell}, X-Y, Z)).$$
Hence
$\roman{e}(A) = \ell$. We furthermore
have the following.

\proclaim{Theorem (5.3)}
Let $Q$ be a parameter ideal in $A$ and
$I = Q : \frak{m}$. Then $\ell_A(I/Q)
\leq 2$. The equality $I^2 = QI$ holds
true if and only if one of the following
conditions is satisfied.
\roster
\item $\ell \geq 2$.
\item $\ell = 1$ and $\ell_A(I/Q) = 1$.
\item $\ell = 1$, $\ell_A(I/Q) = 2$, and
$QB \ne (QB)^{\sharp}$ in $B = A/(x)$.
\endroster
Hence $I^2 = QI$ if either $\ell \geq 2,$
or $\ell = 1$ and $Q \subseteq
\frak{m}^{2}$.
\endproclaim

\demo{Proof}
Let $Q = (f,g)$. Then the sequence $f,g$
is
$B$-regular, so that by (5.1) we get
the exact sequence
$$0 \to A/(\frak{p} + Q) \to A/Q \to
B/QB \to 0. \tag 5.4$$ Hence $\ell_A(I/Q)
\leq 2$, because both the rings
$A/(\frak{p} + Q)
$ and
$B/QB$ are Gorenstein. Since
$A/\frak{p}$ is a DVR and $(Q +
\frak{p})/\frak{p} = (\bar{f},
\bar{g})$, we may assume that 
$(Q + \frak{p})/\frak{p}= (\bar{f}) \ni
\bar{g}$ (here
$\bar{*}$ denotes the reduction mod
$\frak{p}$). Let $\bar{g} =
\bar{c}\bar{f}$ with $c
\in A$. Then, since $Q = (f, g-cf)$,
replacing $g$ by $g - cf$, we get $Q =
(f,g)$ with $g \in \frak{p}$. Since
$\frak{m}/\frak{p} = (\bar{x})$, letting
$\bar{f} =
\bar{\varepsilon}~\bar{x}^n$ with
$\varepsilon \in \roman{U}(A)$ and $n
\geq 1$, we have $Q = (\varepsilon x^n +
a_1, g)$ for some $a_1 \in \frak{p}$. Hence $Q =
(x^n + \varepsilon^{-1}a_1, g)$, so that
$$Q = (x^n + a, b) \tag 5.5$$
with $a, b \in \frak{p}$ and $n \geq 1$.
We then have by (5.4) the exact sequence
$$0 \to A/((x^n) + \frak{p})
\overset{\gamma}\to{\to} A/Q \to B/QB
\to 0, \tag 5.6$$
where $\gamma(1) = x^{\ell}$ mod $Q$. 
We notice that $A/((x^n) +
\frak{p})=R/(X^n, Y, Z)$ is a Gorenstein
ring, containing
$x^{n-1}$ mod $(x^n) + \frak{p}$ as the
non-zero socle. Then by (5.6) 
$\gamma(x^{n-1}~\text{mod}~(x^n) +
\frak{p}) = x^{n+\ell -1}~\text{mod} Q$
is a non-zero element of $I/Q$, 
that is 
$$Q + (x^{n + \ell -1}) \subseteq
I \ \ \  \text{and} \ \ \ x^{n + \ell
- 1}
\not\in Q. \tag 5.7$$
Because $x^{n + \ell -1}a = 0$
(since
$x^{\ell}\frak{p} = (0))$, we get 
$(x^{n + \ell -1})^2 = (x^n +
a)x^{n +
\ell -1}x^{\ell - 1}.$
Hence $(x^{n + \ell -1})^2 \in QI$. This
guarantees that $I^2 = QI$ when
$\ell_A(I/Q) = 1$, because $I =
Q + (x^{n + \ell -1})$ by (5.7).

Now assume that $\ell_A(I/Q) = 2$ and
$\roman{e}(A) = \ell 
\geq 2$. Then $\frak{m}I = \frak{m}Q$ by
Proposition (2.3), whence
$$\mu_A(I) = \ell_A(I/\frak{m}I)
=\ell_A(I/\frak{m}Q) =
\ell_A(I/Q) +\ell_A(Q/\frak{m}Q)
=4, \tag 5.8$$
so that $Q + (x^{n + \ell - 1})
\subsetneq I$. Let $I = Q + (x^{n + \ell
- 1}) + (\xi)$ with $\xi \in A$.
Then, since $B/QB$ is a Gorenstein ring
and the canonical epimorphism $A/Q \to
B/QB$
 in (5.6) is surjective
on the socles, we have $IB = QB + \xi B =
QB : \frak{m}B$. Look at the exact
sequence
$$0 \to A/((x) + Q)
\overset{\delta}\to{\to} B/QB \to
A/((x^{\ell - 1}) + Q) \to 0 \tag 5.9$$
induced from (5.2), where
$\delta(1) = x^{\ell - 1}$ mod $QB$. Then
since
$A/((x) + Q)$ is an Artinian  Gorenstein
ring, choosing
$\Delta \in A$ so that $\frak{m}\Delta
\subseteq (x) + Q$ but $\Delta \not\in
(x) + Q$, by (5.9) we have that $x^{\ell
- 1}\Delta
\not\in QB$ and
$$ IB = QB : \frak{m}B = QB + x^{\ell -
1}\Delta B = QB + \xi B.$$
Let us write $\xi = \varepsilon x^{\ell -
1}\Delta + \rho_0 + x^\ell\varphi_0$ with
$\varepsilon \in \roman{U}(A)$, $\rho_0
\in Q$, and $\varphi_0 \in A$. Then $I =
Q + (x^{n + \ell - 1}) + (\xi)
= Q + (x^{n + \ell - 1}) + (x^{\ell -
1}\Delta + \rho + x^\ell \varphi)$, where
$\rho = \varepsilon^{-1}\rho_0$ and
$\varphi = \varepsilon^{-1} \varphi_0$.
Hence $$I = Q + (x^{n + \ell - 1}) +
(x^{\ell - 1}\Delta + x^{\ell}\varphi)$$
because
$\rho \in Q$. We need the following.

\proclaim{Claim (5.10)}
$\Delta \in \frak{m} = (x) + \frak{p}$.
\endproclaim

\demo{Proof of Claim (5.10)}
Assume $\Delta \not\in \frak{m}$. Then
since $x^{\ell - 1}(\Delta + x \varphi)
\in I$, we have $x^{\ell - 1} \in I$, so
that $I = Q + (x^{\ell - 1})$. This is
impossible, because $\mu_A(I) = 4$ by
(5.8). 
\qed
\enddemo

We write $\Delta = x\sigma + \tau$ with
$\sigma
\in A$ and $\tau \in \frak{p}$. Then
$x^{\ell - 1}\Delta + x^\ell\varphi =
x^{\ell - 1}\tau + x^\ell (\sigma +
\varphi)$ and so 
$$I = Q + (x^{n + \ell - 1}) +
(x^{\ell - 1}\tau + x^{\ell}\varphi_1)\tag
5.11$$
where $\varphi_1 = \sigma + \varphi$.
Suppose that $\varphi_1 \not\in \frak{p}$
and write $\varphi_1 = \varepsilon_1x^q +
\psi_1$ with $\varepsilon_1 \in
\roman{U}(A)$, $q \geq 1$, and $\psi_1
\in \frak{p}$. Then $x^{\ell - 1}\tau +
x^\ell\varphi_1 = x^{\ell - 1}\tau +
\varepsilon_1x^{q + \ell}$ because
$x^\ell \frak{p}= (0)$.
Therefore, letting $\tau_1 =
\varepsilon_1^{-1}\tau$, we get $$I = Q +
(x^{n + \ell - 1}) + (x^{\ell - 1}\tau_1
+ x^{q+\ell}).$$ Because $x^{\ell}\tau_1 =
0$, we have $x^{q+\ell + 1} = x(x^{\ell -
1}\tau_1 + x^{q + \ell} )$, so that $q +
\ell + 1 > n + \ell - 1$ since $\mu_A(I)
= 4$ (otherwise, $I = Q + (x^{\ell - 1}\tau_1
+ x^{q+\ell})$).
Consequently
$x^{q +
\ell}= x^{n +
\ell - 1}(x^{(q+ \ell) - (n + \ell - 1)})$
and so $I = Q + (x^{n + \ell - 1}) +
(x^{\ell - 1}\tau_1)$ with $\tau_1 \in
\frak{p}$. Thus in the expression
(5.11) of $I$ we may assume that
$\varphi_1 \in
\frak{p}$, whence $$I = Q + (x^{n + \ell -
1}) + (x^{\ell - 1}\tau)$$ with
$\tau \in \frak{p}$. Therefore $I^2 = QI
+ (x^{n + \ell - 1}, x^{\ell - 1}\tau)^2 =
QI$, because $(x^{n + \ell - 1})^2 \in QI$
by (5.7) and $x^{\ell - 1}\tau (x^{n +
\ell - 1}, x^{\ell -1 }\tau) = (0)$
(since $x^{\ell} \frak{p} =
(0))$. Thus $I^2 = QI$, if $\ell \geq 2$
or if $\ell = 1$ and $\ell_A(I/Q) = 1$. 

We now consider the case where
$\roman{e}(A) = \ell = 1$ and
$\ell_A(I/Q) = 2$. Our ideal $I$ has in
this case the following normal form
$$I = Q + (x^n, \xi)$$
where $\xi \in \frak{p}$. In fact, $Q +
(x^n) \subseteq I$ and $x^n \not\in Q$ by
(5.7). Since $\ell_A(I/Q) = 2$, the
canonical epimorphism
$A/Q
\to B/QB$ in (5.6) is surjective on the
socles. Hence $IB = QB : \frak{m}B
\supsetneq QB$. Let $I = Q + (x^n) +
(\xi)$ with $\xi
\in A$. If $\xi \not\in \frak{p}$,
letting $\xi = \varepsilon x^q + \xi_1$
with $\varepsilon \in \roman{U}(A)$, $q
\geq 1$, and $\xi_1 \in \frak{p}$, we get
$x\xi = \varepsilon x^{q+1} \in Q$
(recall that $x\frak{p} = (0)$,
since $\ell = 1$). Hence $x^{q+1} \in Q$,
so that $\bar{x}^{q+ 1} \in (\bar{x}^n)
= (Q + \frak{p})/\frak{p}$ in the DVR
$A/\frak{p}$ (cf. (5.5)). Thus $q + 1 \geq
n$. If
$q + 1 = n$, then $x^n \in Q$, which is
impossible by (5.7). Hence $q \geq n$,
and so
$$I = Q + (x^n) + (\varepsilon x^q +
\xi_1) = Q + (x^n, \xi_1)$$ with $\xi_1
\in
\frak{p}$. Thus, replacing $\xi$ by
$\xi_1$ in the case where $\xi \not\in
\frak{p}$, we get
$$I = Q + (x^n, \xi) = (x^n,a,b,\xi)\tag
5.12$$ with $a, b, \xi \in \frak{p}$.
If $QB \ne (QB)^{\sharp}$ in the regular
local ring 
$B = A/(x)$, we have $(IB)^2 =
QB{\cdot}IB$ by Theorem (1.1), since $IB =
QB :
\frak{m}B$. Hence by (5.12) $$(\bar{a},
\bar{b},
\bar{\xi})^2 = (\bar{a},
\bar{b})(\bar{a}, \bar{b},
\bar{\xi})$$
in $B$, where $\bar{*}$ denotes the
reduction mod $(x)$. Therefore
$$(a,b,\xi)^2 \subseteq (a,b)(a,b,\xi) +
(x)$$whence
$$(a,b,\xi)^2 = (a,b)(a,b,\xi)\tag 5.13$$
because $(a,b, \xi) \subseteq
\frak{p}$ and $(x)
\cap
\frak{p}= (0)$. Since $\xi^2 \in
(a,b)(a,b,\xi) = (x^n+a, b)(a,b,\xi)
\subseteq QI$ by (5.13) and
$x^{2n} = (x^n+a)x^n \in QI$, we get that 
$(x^n,\xi)^2 \subseteq QI$, and so $I^2 = QI$ because
$I^2 = QI + (x^n,\xi)^2$ (cf. (5.12)). Thus $I^2 = QI$, if
$QB
\ne (QB)^{\sharp}$. Conversely, assume that
$I^2=QI$. Then $IB \subseteq
(QB)^{\sharp}$, whence $QB \ne (QB)^{\sharp}$
because $QB \subsetneq IB = QB : \frak{m}B
\subseteq (QB)^{\sharp}$. Thus $I^2=QI$
if and only if $QB \ne (QB)^{\sharp}$,
provided $\ell = 1$ and $\ell_A(I/Q) =
2$. This completes the proof of Theorem
(5.3).
\qed
\enddemo

\proclaim{Corollary (5.14)}
Let $\ell = 1$ and $\ell_A(I/Q) = 2$.
Then $I \subseteq Q^{\sharp}$ if and only if
$QB \ne (QB)^{\sharp}$. When this is the
case, the equality $I^2 = QI$ holds true.
\endproclaim

\demo{Proof}
Suppose that $QB = (QB)^{\sharp}$ and $I
\subseteq Q^{\sharp}$. Then $IB = QB$, so that
the monomorphism $A/(\frak{p} + Q) \to
A/Q$ in (5.4) has to be bijective
on the socles, whence $\ell_A(I/Q) = 1$.
This is impossible. If $QB \ne (QB)^{\sharp}$,
we get by Theorem (5.3) that $I^2 = QI$
whence $I \subseteq Q^{\sharp}$.
\qed
\enddemo

Assume that $\ell = 1$ and let
$Q=(x-y,y^2-z^2)$. Then  $\ell_A(I/Q) =
2$. We have by (5.14) $I \not\subseteq
Q^{\sharp}$, since
$QB = (QB)^{\sharp}$ (cf. Theorem
(1.1)). This shows the equality $I^2=QI$
does not necessarily hold true when $\ell
= 1$.

Secondly, let $\frak{a} = (X^3, XY,
Y^2-XZ)$ and let $A = R/\frak{a}$. Let
$x,y$ and $z$ denote the reduction of $X,
Y$ and $Z$ mod $\frak{a}$. Let
$\frak{p} = (x,y)$. We then have the
following.

\proclaim{Lemma (5.15)}
$A$ is a Buchsbaum local ring with $\dim
A = 1$, $\roman{H}^0_{\frak{m}}(A) =
(x^2) \ne (0)$, and $\roman{e}(A) =
\roman{r}(A) = 3.$
\endproclaim

\demo{Proof}
We have $\sqrt{\frak{a}}=(X,Y)$, whence
$\dim A= 1$ and $\roman{Min} A =
\{\frak{p}\}$. We certainly have that
$\frak{m}x^2 = (0)$ and $x^2 \ne 0$. Thus
$(x^2) \subseteq
\roman{H}_{\frak{m}}^0(A)$. Let $$B =
A/(x^2) \cong R/(X^2, XY, Y^2-XZ).$$ We will
show that $B$ is a Cohen-Macaulay ring
with $\roman{e}(B) = 3$. Let $\frak{b} = (X^2, XY, Y^2-XZ)$
and $P = (X,Y)$. Then $P =
\sqrt{\frak{b}}$, $PR_P = (X -
\frac{Y^2}{Z},Y)R_P$, and $\frak{b}R_P =
(X -
\frac{Y^2}{Z},Y^3)R_P$. Hence
$\roman{e}(B) =
\ell_{R_P}(R_P/\frak{b}R_P) = 3$, because
$R/P$ is a DVR. Since
$\frak{n}^2 = Z\frak{n} + \frak{b}$, the
ideal $zB$ is a minimal reduction of the
maximal ideal
$\frak{n}/\frak{b}$
in
$B$, so that we
have
$\roman{e}_{zB}^0(B) = \roman{e}(B) =
3$, while
$\ell_B(B/zB) = \ell_R(R/(X^2, XY, Y^2,
Z))=3$. Thus $\ell_B(B/zB) =
\roman{e}_{zB}^0(B) = 3$, whence
$B=A/(x^2)$ is a Cohen-Macaulay
ring and
$\roman{H}_{\frak{m}}^0(A) = (x^2)$. Let $a \in \frak{m}$ be a
parameter in
$A$. Then $(0) : a
\subseteq
\roman{H}_{\frak{m}}^0(A) = (x^2)$,
since $a$ is a
non-zerodivisor in the
Cohen-Macaulay ring $B =
A/\roman{H}_{\frak{m}}^0(A)$.
Hence $\frak{m}{\cdot}[(0) :
a]=(0)$, so that $A$ is a
Buchsbaum ring. We have
$\mu_{\hat{A}}(\roman{K}_{\hat{A}})
=\mu_{\hat{B}}(\roman{K}_{\hat{B}}) = \roman{r}(B) = 2$,
because
$\roman{H}_{\frak{m}}^1(A)
\cong \roman{H}_{\frak{m}}^1(B)$ and $(X^2, XY, Y^2,
Z) : \frak{n} = \frak{n}$.
Hence $\roman{r}(A) =
\ell_A(\roman{H}_{\frak{m}}^0(A)) +
\roman{r}(B) = 1 + 2 = 3$.
\qed
\enddemo

Let $Q = (a)$ be a parameter ideal in $A$
and put $I = Q : \frak{m}$. Since
$A/\frak{p}$ is a DVR with $z$ mod
$\frak{p}$ a regular parameter, we may
write $a = \varepsilon z^n + b_0$ with
$\varepsilon \in \roman{U}(A)$, $n
\geq 1$, and $b_0 \in \frak{p}$. Hence $Q
= (z^n + b)$, where $b =
\varepsilon^{-1}b_0 \in \frak{p}$.
Consequently, letting $b = xf + yg$ with
$f, g \in A$, we may assume from the
beginning that
$$a = z^n + xf + yg \ \ \ \text{and} \
\ \ Q = (a).\tag 5.16$$

With this notation we have the
following.

\proclaim{Theorem (5.17)}
The equality $I^2 = QI$ holds true if and
only if one of the following conditions
is satisfied.
\roster
\item $f \not\in \frak{m}$.
\item $f \in \frak{m}$ and $n > 1$. 
\endroster
We have $I^3 = QI^2$ but $I^2 \ne QI$, if
$f \in \frak{m}$ and $n = 1$.
\endproclaim

\demo{Proof}
(1) If $f \not\in \frak{m}$, then $A/Q$ is
a Gorenstein ring and $I = Q + (x^2)$. In
fact, choose $F, G \in R$ so that $f, g$
are the reductions of
$F, G$ mod $\frak{a}$, 
respectively. Then $F \not\in \frak{n}$.
We put $V= Z^n + XF + YG$ and $\frak{q} =
(V, XY, Y^2 - XZ)$. Then $\sqrt{\frak{q}} = \frak{n}$ and
so $\frak{q}$ is a parameter ideal in $R$. Let
$x, y$, and
$z$ be, for the moment, the reductions of $X,
Y$, and $Z$ mod $\frak{q}$. We put
$\xi = -F$ mod $\frak{q}$ and $\eta =
G$ mod $\frak{q}$. Then since
$x\xi = z^n + y\eta$, we have
$$
\align
(x\xi)^3 &=(z^n + y\eta)(x\xi)^2\\
&= z^n (x\xi)^2 \ \ \ \ \ (\text{since} \
\  xy=0)\\ &=(x\xi{\cdot}z)(x\xi)z^{n-1}\\
&=(y^2\xi)(x\xi)z^{n-1} \ \ \
\ \ (\text{since}
\
\  y^2 = xz)\\ &=0.
\endalign
$$
Thus $x^3 = 0$ in $R/\frak{q}$.
Consequently $X^3 \in \frak{q}$, so that
$\frak{q}=(V, X^3, XY, Y^2-XZ)$. Hence
$A/Q = A/(z^n + xf + yg) 
\cong  R/(V, X^3, XY, Y^2-XZ)=
R/\frak{q}$ and so $A/Q$ is
a Gorenstein ring. Since $\ell_A(I/Q) = 1$
and $x^2
\not\in Q$ (otherwise, $x^2 \in
\roman{H}_{\frak{m}}^0(A) \cap Q =
(0)$; recall that $A$ is a Buchsbaum
ring), we get that $I = Q + (x^2)$. Thus
$I^2 = QI$.

(2) Suppose that $f \not\in \frak{m}$ and
$n > 1$. Then, since $xa = xz^n$ and $ya
= yz^n + y^2g = yz^n + xzg$, we get 
$$a\frak{p} =(xz^n, yz^n + y^2g)
\subseteq (z)\tag 5.18$$
and $\frak{m}{\cdot}(xz^{n-1}, x^2)
\subseteq a\frak{p}$. We claim that the
reductions of $xz^{n-1}$ and $x^2$ mod
$a\frak{p}$ are linearly independent in
$\frak{p}/a\frak{p}$ over the field
$A/\frak{m}$. In fact, let $c_1, c_2 \in
A$ and assume that $c_1(xz^{n-1}) + c_2x^2
\in a\frak{p}$. Then since $n > 1$ and $a\frak{p} \subseteq
(z)$ by (5.18), we have $c_2x^2 \in (z)$,
and so
$c_2x^2
\in
\roman{H}_{\frak{m}}^0(A) \cap (z) = (0)$
(recall that $(z)$ is a parameter ideal
in $A$). Hence $c_2 \in \frak{m}$ so that 
$c_1(xz^{n-1}) \in a\frak{p}$.
Suppose $c_1 \not\in \frak{m}$ and write
$xz^{n-1} = xz^n\varphi + (yz^n +
y^2g)\psi$ with $\varphi, \psi \in A$.
Then because $xz^{n-1}(1 - z \varphi) =
(yz^n + y^2g)\psi$, we get $xz^{n-1} =
(yz^n + y^2g)\rho$ for some $\rho \in A$.
Hence
$$z^{n-1}(x-yz\rho) = y^2g\rho= xzg\rho.\tag
5.19$$
Now notice that $A/(x) \cong R/(X,Y^2)$
and we see that $z$ is $A/(x)$-regular.
Because $z^{n-1}(-yz\rho) \equiv 0$ mod
$(x)$ (cf. (5.19)), we get $y\rho \equiv
0$ mod
$(x)$, whence $y^2\rho = 0$. This implies
by (5.19) that $$x -yz\rho \in (0) :
z^{n-1} = (0) : z = (x^2)$$ since $z$ is a parameter in our
Buchsbaum ring $A$. Thus $x \in \frak{m}^2$
which is impossible. Hence $c_1 \in
\frak{m}$.

Now let $B = A/\frak{p}$ and look at the 
canonical exact sequence
$$0 \to \frak{p}/a\frak{p} \to A/Q \to
B/QB \to 0\tag 5.20$$
of $A$-modules and we have
$$2 \leq \ell_A((0) :_{\frak{p}/a\frak{p}}
\frak{m}) \leq \ell_A(I/Q) \leq
\roman{r}(A) = 3.\tag 5.21$$ If
$\ell_A(I/Q) =
\roman{r}(A) = 3$, then $I^2 = QI$ by
Theorem (3.9). Hence to prove $I^2 = QI$,
we may assume $\ell_A(I/Q) \leq 2$.
Therefore $\ell_A((0)
:_{\frak{p}/a\frak{p}}
\frak{m}) = \ell_A(I/Q) = 2$  by
(5.21) so that by (5.20) we have 
$I = Q + (xz^{n-1},x^2)$,
because $[(0)
:_{\frak{p}/a\frak{p}}
\frak{m}]$ is
generated by the
reductions of $xz^{n-1}$ and $x^2$ mod
$a\frak{p}$. Hence $I^2 = QI +
(xz^{n-1},x^2)^2 = QI$, since
$x^2\frak{m} = (0)$.

(3) Suppose that $f \in \frak{m}$ and $n =
1$. Let $f = xf_1+yf_2+zf_3$ with $f_i
\in A$. Then $a = z + xf + yg = z +
x^2f_1 + y(g + yf_3)$, because $y^2 = xz$.
Consequently, replacing $f$ by
$xf_1$ and $g$ by $g + yf_3$, we may
assume in the expression (5.16) of $I$
that
$$a = z +  x^2f + yg \ \ \ \text{and} \ \
\  Q = (a).$$
Hence $a\frak{p} = (xz, yz + y^2g) = (xz,
yz) =z\frak{p}$ (recall that $y^2 = xz$).
Look at the exact sequence 
$$0 \to \frak{p}/a\frak{p} \to A/(z)
\to B/zB \to 0\tag 5.22$$
of $A$-modules. Then, because $A/(z)
\cong R/(X^3, XY, Y^2, Z)$, we see
$\ell_A(((z) : \frak{m})/(z)) = 2$ and
$(z) : \frak{m} = (z) + (x^2, y)
\subseteq (z) + \frak{p}$. Hence
in (5.22) the canonical
epimorphism
$A/(z) \to B/zB$ is zero on the socles.
Thus $\ell_A((0) :_{\frak{p}/a\frak{p}}
\frak{m}) = 2$ and $[(0)
:_{\frak{p}/a\frak{p}} \frak{m}]$ is
generated by the reductions of $x^2$ and
$y$ mod $a\frak{p} = z\frak{p}$.
Consequently $Q + (x^2,y)
\subseteq I$ by (5.20).

\proclaim{Claim (5.23)}
$\ell_A(I/Q) \ne 3$.
\endproclaim

\demo{Proof of Claim (5.23)}
Assume $\ell_A(I/Q) = 3$. Then $I^2
= QI$ by Theorem (3.9), since
$\ell_A(I/Q) =
\roman{r}(A)$. Thus
$IB = QB$, because
$IB \subseteq (QB)^{\sharp} = QB$ (notice that
$B$ is a DVR). Hence in (5.20) the
epimorphism $A/Q \to B/QB$ has to be zero
on the socles, and so
$\ell_A(I/Q) = \ell_A((0)
:_{\frak{p}/a\frak{p}}
\frak{m}) = 2$,
which is impossible.
\qed
\enddemo

By this claim we see that $I = Q + (x^2,
y)$, whence $I^2 = QI + (y^2)$.
Consequently, $I^3 = QI^2$, because $y^3
= y{\cdot}xz = 0$. In contrast, $I^2 \ne
QI$, because $y^2 \not\in QI$. To see
this, assume that $y^2 \in QI$ and choose
$F, G \in R$ so that $f, g$ are the
reductions of $F,G$ mod $\frak{a}$, respectively. Let
$K=(Z^2+YZG, YZ+Y^2G, X^3,XY,Y^2-XZ)$.
Then $Y^2 \in K$, because
$QI = (z+x^2f+yg)(z,x^2,y)
=(z^2+yzg, yz+y^2g).$
Hence $$K = (X^3, Y^2,
Z^2, XY, YZ, ZX)$$ which is impossible,
since
$\mu_R((X^3, Y^2,
Z^2, XY, YZ, ZX)) = 6$ while
$\mu_R(K) \leq 5$. Thus $y^2 \not\in
QI$, which completes
the proof of Theorem (5.17).
\qed
\enddemo

If $Q \subseteq \frak{m}^2$, then $n \geq
2$, and so by Theorem (5.17) we
readily get the following.

\proclaim{Corollary (5.24)}
$I^2 =QI$ if $Q \subseteq \frak{m}^2$.
\endproclaim



\Refs

\widestnumber\key{CHV}

\ref
\key C
\by N.~T.~Cuong
\paper P-standard systems of parameters and p-standard ideals in local
rings
\yr 1995
\pages 145-161
\vol  20
\jour Acta. Math. Vietnamica
\endref

\ref
\key CHV
\by A. Corso, C. Huneke, and W. V. Vasconcelos
\paper On the integral closure of ideals
\yr 1998
\pages 331-347
\vol  95
\jour manuscripta math.
\endref

\ref
\key CP
\by A. Corso and C. Polini
\paper Links of prime ideals and their Rees algebras
\yr 1995
\pages 224-238
\vol  178
\jour J. Alg.
\endref



\ref
\key CPV
\by A. Corso,  C. Polini, and W. V. Vasconcelos
\paper Links of prime ideals
\yr 1994
\pages 431-436
\vol  115
\jour Math. Proc. Camb. Phil. Soc.
\endref

\ref
\key CST
\by     N. T. Cuong, P. Schenzel and N. V. Trung
\paper  Verallgemeinerte Cohen-Macaulay-Moduln 
\jour   M. Nachr.
\vol    85
\yr     1978
\pages  57-73
\endref


\ref
\key    G1
\by     S.~Goto
\paper  On Buchsbaum rings
\jour   J. Alg.
\vol    67
\yr     1980
\pages  272-279
\endref

\ref
\key    G2
\by     S. Goto
\paper  On the associated graded rings of parameter ideals 
        in Buchsbaum rings
\jour   J. Algebra
\vol    85
\yr     1983
\pages  490-534
\endref

\ref
\key G3
\by S.~Goto
\paper Integral closedness of complete-intersection ideals
\yr 1987
\pages 151-160
\vol  108
\jour J. Alg.
\endref

\ref
\key GH
\by S. Goto and F. Hayasaka
\paper Finite homological dimension and
primes associated to integrally closed
ideals II 
\jour J. Math. Kyoto Univ. 
\toappear
\endref

\ref
\key GN
\by S.~Goto and K.~Nishida
\paper Hilbert coefficients and
Buchsbaumness of associated graded rings
\yr 2003
\pages 61-74
\vol  181
\jour J. Pure and Applied Alg.
\endref



\ref
\key GSa
\by S.~Goto and H.~Sakurai
\paper The equality $I^2=QI$ in Buchsbaum
 rings with multiplicity 2
\jour Preprint 2002 
\endref

\ref
\key GSu
\by S.~Goto and N.~Suzuki
\paper Index of reducibility of
parameter ideals in a local ring
\yr 1984
\pages 53-88
\vol  87
\jour J. Alg.
\endref

\ref
\key H
\by C.~Huneke
\paper The theory of $d$-sequences and
powers of ideals
\yr 1982
\pages 249-279
\vol  46
\jour Ad. in Math.
\endref


\ref
\key    HK
\by    J. Herzog and E. Kunz (eds.)
\book Der kanonische Modul eines Cohen-Macaulay-Rings
\bookinfo Lecture Notes in Math.
\vol 238
\publ Springer-Verlag
\publaddr Berlin ${\cdot}$  Heidelberg  ${\cdot}$   New York 
${\cdot}$  Tokyo
\yr 1971
\endref


\ref
\key    SV1
\by     J. St\"{u}ckrad and W. Vogel
\paper  Eine Verallgemeinerung der Cohen-Macaulay-Ringe und
Anwendungen auf ein Problem der Multiplizit\"{a}tstheorie
\jour   J. Math. Kyouto Univ.
\vol    13
\yr     1973
\pages  513-528
\endref


\ref
\key    SV2
\by     J. St\"{u}ckrad and W. Vogel
\book   Buchsbaum rings and applications
\publ Springer-Verlag
\publaddr  Berlin, New York, Tokyo
\yr 1986
\endref

\ref
\key    Y1
\by     K. Yamagishi
\paper  The associated graded modules of 
Buchsbaum modules with respect to 
$\frak{m}$-primary ideals in the
equi-$\Bbb{I}$-invariant case
\jour   J. Alg.
\vol    225
\yr     2000
\pages  1-27
\endref

\ref
\key    Y2
\by     K. Yamagishi
\paper  Buchsbaumness in Rees modules associated to ideals of
minimal multiplicity in the equi-$\Bbb
I$-invariant case
\jour   J. Alg.
\vol    251
\yr     2002
\pages  213-255
\endref




\endRefs
\enddocument